\documentclass[times,final,nopreprintline]{elsarticle}

\usepackage{framed,multirow}
\usepackage{lineno}
\usepackage{amssymb,amsthm}
\usepackage{latexsym}
\usepackage{bbm}
\usepackage{ulem}
\usepackage{dsfont}
\usepackage{hyperref}
\usepackage{amsfonts}
\usepackage{graphicx}
\usepackage{epstopdf}
\usepackage{overpic}
\usepackage{amsmath,esint}
\usepackage{epsfig}
\usepackage{todonotes} 
\usepackage{tikz}
\usepackage{cleveref}
\usetikzlibrary{arrows}
\usetikzlibrary{fadings}
\usepackage{pgfplots}
\pgfplotsset{compat=1.10}
\usepgfplotslibrary{fillbetween}
\usetikzlibrary{patterns}
\usepackage{algorithm}
\usepackage{appendix}
\newtheorem{remark}{Remark}

\usepackage[utf8]{inputenc}
\usepackage{bbm}
\usepackage{enumitem}
\usepackage{mathrsfs}
\usepackage{verbatim}
\usepackage{graphicx}
\usepackage{algorithm}
\usepackage{algorithmicx}
\usepackage{algpseudocode}
\usepackage{mathtools}
\usepackage{pifont}
\usepackage{calc}
\usepackage{accents}
\usepackage{cancel}
\usepackage[margin=1in,footskip=0.25in]{geometry}

\setenumerate{label=(\roman*),itemsep=3pt,topsep=3pt}

\newcommand{\RN}[1]{%
  \textup{\uppercase\expandafter{\romannumeral#1}}%
}

\DeclarePairedDelimiter{\ceil}{\lceil}{\rceil}
\setenumerate{label=(\roman*),itemsep=3pt,topsep=3pt}

\usepackage{ifpdf}
\usepackage{array}

\newcommand{\bbR}{\mathbb R}

\newcommand{\bgamma}{\boldsymbol \gamma}

\newcommand{\bGamma}{\mbox{\boldmath{$\Gamma$}}}

 \newcommand{\bd}{\mathbf d}

 \newcommand{\cJ}{\mathcal J}
 
\newcommand{\cM}{\mathcal M} \newcommand{\cN}{\mathcal N}
\newcommand{\cO}{\mathcal O}  
 
\newcommand{\cS}{\mathcal S}

\newcommand{\Rthree}{{\bbR^3}}
\newcommand{\Stwo}{{\cS^2}}

\DeclareMathOperator*{\argmin}{arg\,min}

\graphicspath{{Fig/}}

\journal{Journal of Computational Physics}

\begin{document}

\begin{frontmatter}

\title{Adjoint DSMC for nonlinear Boltzmann equation constrained optimization}


\author[address1]{Russel  Caflisch}
\ead{caflisch@courant.nyu.edu}
\address[address1]{Courant Institute of Mathematical Sciences, New York University, New York, NY 10012.}
\author[address1]{Denis Silantyev}
\ead{silantyev@courant.nyu.edu}
\author[address1]{Yunan Yang\corref{mycorrespondingauthor}}
\ead{yunan.yang@nyu.edu}
\cortext[mycorrespondingauthor]{Corresponding author}

\begin{abstract}
Applications for kinetic equations such as optimal design and inverse problems often involve finding unknown parameters through gradient-based optimization algorithms. Based on the adjoint-state method, we derive two different frameworks for approximating the gradient of an objective functional constrained by the nonlinear Boltzmann equation. While the forward problem can be solved by the DSMC method, it is difficult to efficiently solve the high-dimensional continuous adjoint equation obtained by the ``optimize-then-discretize'' approach. This challenge motivates us to propose an adjoint DSMC method following the ``discretize-then-optimize'' approach for Boltzmann-constrained optimization. We also analyze the properties of the two frameworks and their connections. Several numerical examples are presented to demonstrate their accuracy and efficiency.
\end{abstract}



\end{frontmatter}




\section{Introduction}
The development of modern technology requires accuracy in modeling physical processes. One critical task is to model the kinetic behavior of rarefied gas, which is required for low-pressure gas flow and cannot be accurately described by the Navier--Stokes equations. The Boltzmann equation models the dynamics of a many-particle system through a velocity distribution function, with a nonlinear collision operator that describes binary particle interactions. The Boltzmann equation is a powerful tool from the kinetic theory describing molecular gas dynamics, radiative transfer, plasma physics, and grain and polymer flow~\cite{alonso2015boltzmann}.

The Boltzmann equation can also be used for design, optimization, control, and inverse problems. A few of the many examples include the design of the semiconductor device, the topology optimization of the gas flow channel, and risk management in quantitative finance~\cite{choulli1996inverse,dragulescu2003applications,cheng2011recovering,sato2019topology,lai2020reconstruction}. One may use a kinetic perspective to tackle optimal control problems for a large system of interacting agents. Common kinetic models involve the Boltzmann model and some mean-field models~\cite{toscani2006kinetic,caponigro2013sparse,albi2014boltzmann,fornasier2014mean2,fornasier2014mean,albi2015kinetic}. Many of these applications involve finding unknown or optimal parameters in  Boltzmann-type equations such that an objective function formed by the computational or experimental data is optimized, i.e., optimization problems with  PDE constraints. However, due to the challenges caused by the complex nonlinear collision term, the linear Boltzmann equation~\cite{choulli1996inverse} or simplified collision operators such as the Bhatnagar–Gross–Krook (BGK) model~\cite{bhatnagar1954model,sato2019topology}, the Ellipsoidal Statistical model~\cite{holway1966kinetic}, or the Shakhov model~\cite{shakhov1968generalization} have been used as alternatives. Although these simplified models provide good approximations for many scenarios, the original nonlinear collision operator is still preferred for better accuracy~\cite{cercignani1988boltzmann}.

\begin{figure}
\begin{center}
{\footnotesize{
\tikzstyle{decision1} = [rectangle, draw, fill=blue!20,
    text width=6em, text centered, node distance=3cm, rounded corners, minimum height=4em]
\tikzstyle{decision2} = [rectangle, draw, fill=red!20,
    text width=6em, text centered, node distance=3cm, rounded corners, minimum height=4em]
\tikzstyle{block} = [rectangle, draw, fill=orange!20,
    text width=5em, text centered, node distance=3cm, rounded corners, minimum height=4em]
\tikzstyle{block1} = [rectangle, draw, fill=white!20,
    text width=5em, text centered, node distance=3cm, rounded corners, minimum height=4em]
\tikzstyle{line} = [draw, thick, color=black, -latex']
\tikzstyle{line2} = [draw, thick,dash dot, color=black, -latex']
\tikzstyle{line3} = [draw, thick,dash dot, color=red, -latex']

\hspace*{18pt}\begin{tikzpicture}[scale=2, node distance = 1cm, auto]
    \node(OTD1) [block,text width=10em,node distance=1cm]  {the Boltzmann-constrained optimization problem~\eqref{eq:OTD_obj}} ;
    \node (OTD2) [decision1, above of=OTD1,text width=9em,node distance=2cm, xshift = 2cm]  {the continuous optimality condition~\eqref{optimization1} and the adjoint equation~\eqref{backwards1}};
    \node (DTO1) [decision2, below of=OTD1,xshift = 2cm, node distance=2cm, text width=9em]  {a Lagrangian~\eqref{eq:DTO_obj} based on the DSMC method for the state equation~\eqref{eq:homoBoltz}};
    \node (OTD3) [decision1, right of=OTD2,node distance=4.5cm,text width=9em] {solve state and adjoint equations~\eqref{eq:homoBoltz} \&~\eqref{backwards1} numerically and compute~\eqref{optimization1}};
    \node (DTO2) [decision2, right of=DTO1, node distance=4.5cm,text width=9em] {the discrete optimality~\eqref{DalphaEqtn} and the adjoint DSMC equation~\eqref{eq:AB_gamma}};
        \node (grad) [block, above of=DTO2, xshift = 2cm, node distance=2cm,text width=9em] {the gradient of~\eqref{eq:OTD_obj}};
\path [line] (OTD1) |- node [anchor=east,yshift = -0.7cm]{
\begin{tabular}{cc}
optimize \\
OTD\\
\end{tabular}
}(OTD2);
\path [line] (OTD2) -- node{discretize}(OTD3);
\path [line] (OTD3) -| node[anchor=west,yshift = -0.7cm]{approximate}(grad);
\path [line2] (OTD1) |-  node [anchor=east,yshift = 0.7cm]{
\begin{tabular}{cc}
discretize \\
DTO\\
\end{tabular}
}(DTO1);
\path [line2] (DTO1) -- node{optimize}(DTO2);
\path [line2] (DTO2) -| node[anchor=west,yshift = 0.7cm]{approximate}(grad);
\end{tikzpicture}
}}
\end{center}
\caption{A diagram summarizing the OTD approach (solid line) and the DTO approach (dash line) to compute the gradient with respect to the unknown parameter for the optimization problem~\eqref{eq:OTD_obj} constrained by the Boltzmann equation~\eqref{eq:homoBoltz}.}\label{fig:flowchart}
\end{figure}
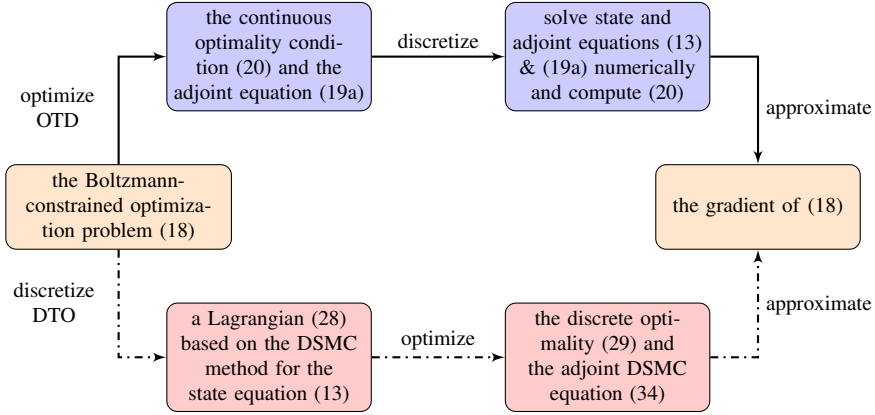

The main contribution of this work involves the derivation of two optimization frameworks based on the spatially homogeneous Boltzmann equation with the nonlinear collision operator and the design of an efficient adjoint DSMC method for the gradient calculation. In this paper, we focus on the collision kernel for Maxwell molecules, but we expect our method to be generalized to other collision kernels. The frameworks should be broadly applicable to optimal control, optimal design, and general computational inverse problems. We employ two different approaches to derive numerical algorithms: optimize-then-discretize (OTD) and discretize-then-optimize (DTO)~\cite{hinze2012discretization}. Fig.~\ref{fig:flowchart} summarizes these two approaches. In the OTD approach, we obtain a continuous adjoint equation as an optimality condition, which we discretize by a grid-based method or a Monte Carlo type method. In the DTO approach, we discretize the Boltzmann equation and the objective function using the direct simulation Monte Carlo method (DSMC)~\cite{bird1970direct,nanbu1980direct,caflisch1998monte}, from which a discrete optimality condition and the adjoint DSMC system are derived. The adjoint DSMC from the DTO approach is radically more cost-effective than the continuous adjoint equation in the OTD approach. To our knowledge, this is the first efficient numerical scheme to compute the gradient for nonlinear Boltzmann equation constrained optimization. In addition to the derivations, we also analyze the properties of the two adjoint systems and investigate connections between them. The analysis provides a better understanding of the adjoint-state method for constrained optimization problems. 

The paper is arranged as follows. In~\Cref{sec:background}, we briefly review the homogeneous Boltzmann equation with a nonlinear collision operator, the linearization of the collision operator, and the classical DSMC algorithm for the numerical solution. We derive the continuous adjoint equation for the continuous Boltzmann equation in~\Cref{ContinuousFormulation} following the Lagrangian method. A particle method is then proposed to solve the continuous adjoint equation numerically.
In~\Cref{DiscreteFormulation}, we regard the DSMC discretization of the Boltzmann equation as the forward problem and the particle velocities as the state variables. We then obtain equations for the adjoint particles following the adjoint-state method. Since one can compute the gradient of the optimization problem by both the continuous and the discrete adjoint systems, we analyze their properties and prove the critical connections between the two adjoint systems in~\Cref{sec:discussion}. Finally, we show several numerical examples in~\Cref{sec:numerics} to demonstrate the accuracy of the gradient formulae obtained by both the continuous and the DSMC adjoint systems. Detailed comparison in terms of memory, accuracy, and computational cost of different numerical schemes are presented.
Conclusion and future research directions follow in~\Cref{sec:conclusion}. Two appendices provide additional details on the numerical methods and numerical analysis of the results.



\section{Boltzmann equation and the direct simulation Monte Carlo (DSMC) method}\label{sec:background}
In this section, we first give a short introduction to the Boltzmann equation and some of its relevant properties. The later part of the section is devoted to the description of the classical DSMC method.
\subsection{Boltzmann equation}
We consider the Boltzmann equation
\begin{equation*}
\frac{\partial f}{\partial t}   + v \cdot \nabla_x f =  Q(f,f)
\end{equation*}
with the initial condition 
\begin{equation*}
    f(x,v,t=0) = f_0(x,v),
\end{equation*}
where $f(x,v,t)$ is a nonnegative probability density function that describes the time evolution of the distribution of particles which move with velocity $v\in \bbR^3$ at the position $x \in \bbR^3$ at time $t>0$. The bilinear (nonlinear) collision operator $Q(f,f)$ that describes the binary collisions among particles is defined as follows:
\begin{equation}~\label{eq:nonlinearCollision}
    Q(f,f) = \int_{\bbR^3} \int_{\cS^2} q(v-v_1,\sigma) (f(v_1')f(v') - f(v_1) f(v)) d\sigma dv_1,
\end{equation}
in which $(v',v_1')$ represent the post-collisional velocities associated with the pre-collisional velocities $(v,v_1)$ and the $\sigma$ integral is over the surface of unit sphere $\cS^2$. By conserving the momentum $v+v_1$ and the energy $v^2 + v_1^2$, we have
\begin{subequations}~\label{BoltzmannSolution}
\begin{align}
v' &= 1/2 (v+v_1) + 1/2 |v-v_1| \sigma,  \\
v_1'&=1/2 (v+v_1) - 1/2  |v-v_1| \sigma,
\end{align}
\end{subequations}
where $\sigma$ is a collision parameter representing a unit direction of the relative velocity of particles after collision~\eqref{BoltzmannSolution}. We will hereafter use the shorthand notation $f, f_1, f', f_1'$ to denote $f(v), f(v_1), f(v')$ and $f(v_1')$.

The collision operator $Q$ has the following symmetries, which are related to the conservation of mass, momentum and energy: The transformations 
\begin{align*}
v &\leftrightarrow v_1  \mbox{ and } v' \leftrightarrow v_1'  \\
v &\leftrightarrow v'  \mbox{ and } v_1 \leftrightarrow v_1' 
\end{align*}
with $v'$ and $v_1'$ defined by~\eqref{BoltzmannSolution} are isometries for $(\sigma, v,v_1) \in \cS^2 \times \bbR^3 \times \bbR^3$~\cite{cercignani1969mathematical,villani2002review}. We define the unit vector along the direction of $v-v_1$ as $\alpha = \frac{v-v_1}{|v-v_1|}$. The scattering angle $\theta = \cos^{-1}{(\sigma\cdot \alpha)}$. Also, $|v-v_1|=|v'-v_1'|$, and thus we have
\begin{equation*}
q(v-v_1,\sigma) = \tilde q(|v-v_1|,\theta)=\tilde q(|v'-v_1'|, \theta) = q(v'-v_1',\alpha),
\end{equation*}
where $\tilde{q}$ is related to $q$ by a change of variable~\cite{villani2002review}.
It follows that for any measurable function $F(v,v_1,v',v_1')$ under suitable integrability conditions 
\begin{subequations}\label{eq:symmetry}
\begin{align}
\iiint F(v,v_1,v',v_1')q d\sigma dv dv_1 
&= \iiint F(v_1,v,v_1',v') q d\sigma dv dv_1 \\ 
 &= \iiint F(v',v_1',v,v_1) q d\sigma dv dv_1  \\ 
 &= \iiint F(v_1', v', v_1, v)q d\sigma dv dv_1. 
\end{align}
\end{subequations}
These symmetries are used repeatedly through the paper.

The kernel $q$ is a nonnegative function that characterizes the details of the binary interactions. 
For the Variable Hard Sphere (VHS) model, the collision kernel is
\begin{equation}
    q(v-v_1,\sigma) =  \tilde q(|v-v_1|,\theta) = C_\beta(\theta) |v-v_1|^\beta,
\end{equation}
In particular, when $\beta = 0$, the collision kernel corresponds to the Maxwellian gas, which (along with  $C_\beta$ being constant) is the model that we focus on in this paper. For the Coulomb interaction, the collision kernel is given by the Rutherford formula
$$
q(v-v_1,\sigma) = \tilde q(|v-v_1|,\theta) = \dfrac{1}{|v-v_1|^3\sin^4(\theta/2)}.
$$


\subsection{Operator formulation for collisions}\label{sec:operator_notation}

In this section, we rewrite the collision rules~\eqref{BoltzmannSolution}, in an operator formulation that clarifies some properties that will be useful for the analysis in later sections. First we denote $\alpha = (v -v_1){\hat~}$ and  $\sigma = (v' -v_1'){\hat~}$ using the notation $x{\hat~} = x/|x|$, in which the equation for $\sigma$ and the relation $|v-v_1| = |v'-v'_1|$ follow from~\eqref{BoltzmannSolution}. As described in~\cite[P.53]{villani2002review}, the change of variables   
\begin{equation}
    (v,v_1,\sigma) \longrightarrow     (v',v'_1,\alpha) 
\end{equation}
is an involution  with unit Jacobian.

We can rewrite~\eqref{BoltzmannSolution}, in terms of operators, as
\begin{equation}\label{eq:AB_vel_General}
        \begin{pmatrix}
    v'\\
   v_1'
    \end{pmatrix}    = A(\sigma,\alpha) \begin{pmatrix}
    v\\
    v_1
    \end{pmatrix}, \quad 
   \begin{pmatrix}
    v\\
     v_1
    \end{pmatrix}     = B(\sigma,\alpha)    \begin{pmatrix}
    v'\\
    v_1'
    \end{pmatrix},
\end{equation}
where
\begin{equation} \label{eq:AB_def}
A(\sigma ,\alpha )= \dfrac{1}{2}\begin{pmatrix}
    I+\sigma  \alpha ^T & I-\sigma  \alpha ^T\\
    I-\sigma  \alpha ^T & I+\sigma  \alpha ^T
    \end{pmatrix},\,B(\sigma ,\alpha )= \dfrac{1}{2}\begin{pmatrix}
    I+\alpha  \sigma ^T & I-\alpha  \sigma ^T\\
    I-\alpha  \sigma ^T & I+\alpha  \sigma ^T
    \end{pmatrix},
\end{equation}
where $I$ is the identity matrix in $\Rthree$ and $B= A^T = A^{-1}$.

Now consider small perturbations $\delta v$ and $\delta v_1$ in the pre-collision velocities $v$ and $v_1$. Since the collision parameter vector $\sigma$ is chosen independently of $v$ and $v_1$ for a Maxwellian gas, there is no need to perturb $\sigma$; i.e., $\delta\sigma = 0$.  Also note that
$\alpha = \partial_v |v - v_1|  = - \partial_{v_1} |v - v_1| $. 
The resulting first-order variations $\delta v'$ and $\delta v_1'$ in the post-collision velocities are  
\begin{equation}\label{eq:fwd_perturb_General}
        \begin{pmatrix}
    \delta v'\\
    \delta v_1'
    \end{pmatrix}   =
    A(\sigma,\alpha)    
 \begin{pmatrix}
    \delta  v\\
    \delta v_1
    \end{pmatrix}.
\end{equation}



\subsection{The linearized collision operator}
The linearized collision operator is defined through perturbation theory~\cite{cercignani1969mathematical}, most frequently by linearization around a Maxwellian equilibrium. In this subsection, we define the linearized operator, due to multiplicative perturbation~\cite{arkeryd1988stability,caflisch1980boltzmann} around a non-equilibrium distribution and its adjoint operator under a weighted $L^2$ inner product. 

For a general velocity distribution $f$, consider the multiplicative perturbation $\tilde f = f(1+\psi)$. Since the collision operator $Q(\tilde f,\tilde f)$ defined in~\eqref{eq:nonlinearCollision} is bilinear, 
\begin{align}
    Q(\tilde f,\tilde f)  - Q(f,f)  &=  Q(f, f\psi) + Q(f\psi,f) + Q(f\psi,f\psi) \nonumber
    \\
   &\approx Q(f, f\psi) + Q(f\psi,f), \nonumber
\end{align}
in which $\approx$ means that quadratic terms in $\psi$ are neglected.
We  define the linearized operator $L[f]$ based on $f$ and applied to $\psi$ as
\begin{equation}~\label{eq:linearQ}
L[f]  \psi =  f^{-1} [ Q(f, f\psi) + Q(f\psi,f)].
\end{equation}

In an original paper of Maxwell in 1866~\cite{maxwell1867iv}, the Boltzmann equation was written in the weak formulation~\cite{villani2002review}. Following this approach for a test function $\gamma(v,t)$ and using the symmetries~\eqref{eq:symmetry}, one gets the identity
\begin{equation}~\label{eq:variation0}
    \int_\Rthree Q(f,f)\gamma(v) dv = \frac{1}{2} \int_\Rthree\int_\Rthree\int_\Stwo ff_1 (\gamma' + \gamma_1'-\gamma - \gamma_1)qd\sigma dv_1 dv
\end{equation}
where $\gamma_1,\gamma'$ and $\gamma_1'$ are shorthand notations for $\gamma(v_1,t)$, $\gamma(v',t)$ and $\gamma(v_1',t)$.
We will use (\ref{eq:variation0}) to find the adjoint of $L[f]$.
We  replace $f$ by $f(1+\psi)$ in (\ref{eq:variation0}) and calculate the first-order terms  with respect to $\psi$. As in the derivation of~\eqref{eq:linearQ},  the first-order terms in the left-hand side of~\eqref{eq:variation0} are derived as 
\begin{equation*}
  \int_\Rthree  [Q(f+f\psi,f+f\psi) - Q(f,f)]\gamma(v) dv \approx \int_\Rthree f(L[f] \psi)  \gamma(v) dv.
\end{equation*}   
Similarly, the first-order terms in the right-hand side of~\eqref{eq:variation0} are derived as
\begin{align*}
 \frac{1}{2} \iiint [(f+f\psi)(f_1+f_1\psi_1) - ff_1](\gamma' + \gamma_1'-\gamma - \gamma_1)qd\sigma dv_1 dv  \approx   \iiint ff_1 \psi(\gamma' + \gamma_1'-\gamma - \gamma_1)qd\sigma dv_1 dv.
\end{align*}
Combining both sides, we obtain the first-order variation equation
\begin{equation*}
    \int_\Rthree f \left((L[f] \psi) \gamma  - \psi \int_\Rthree\int_\Stwo f_1 (\gamma' + \gamma_1'-\gamma - \gamma_1)qd\sigma dv_1\right) dv=0,
\end{equation*}
which is equivalent to the following
\begin{equation*} 
    (L[f]\psi, \gamma)_{L^2(\mu)} = (\psi, L^*[f]\gamma)_{L^2(\mu)} ,
\end{equation*}
in which 
\begin{equation}~\label{eq:adj_linear_operator}
L^*[f]\gamma = \int_\Rthree\int_\Stwo f_1 (\gamma' + \gamma_1'-\gamma - \gamma_1)qd\sigma dv_1.
\end{equation}
and $(h_1,h_2)_{L^2(\mu)}  = \int h_1h_2fdv$ is the inner product for the weighted Hilbert space $L^2(\mu)$ in which $f dv = d\mu$. 
This shows that operator $L^*[f]$ 
is the adjoint of  $L[f]$.

In the special case that $f=\mathcal{M}$ is a Maxwellian equilibrium distribution, using the symmetries~\eqref{eq:symmetry} and the equilibrium property $\cM\cM_1 = \cM'\cM_1'$,
\begin{equation}
L[\cM]  \psi = \int_\Rthree \int_\Stwo (\psi' + \psi_1' -\psi-\psi_1 )\cM(v_1)q(v-v_1,\sigma) d\sigma dv_1. \label{eq:linearEqui}
\end{equation}
Comparison of (\ref{eq:linearEqui}) and (\ref{eq:adj_linear_operator}) shows that $L^*[\cM]=L[\cM]$; i.e., that $L[\cM]$ is self adjoint.
The linearized collision operator around a Maxwellian has been extensively studied in the literature~\cite{cercignani1969mathematical,villani2002review,bobylev2008some}.




\subsection{The direct simulation Monte Carlo (DSMC) method} 
In this section, we describe the classical DSMC method~\cite{bird1970direct} for the spatially homogeneous Boltzmann equation following the presentation of~\cite{pareschi2001introduction}:
\begin{equation}\label{eq:homoBoltz}
    \frac{\partial f}{\partial t}  =  Q(f,f).
\end{equation}
As stated above, our focus is also on a Maxwellian gas for which the collision kernel does not depend on the relative velocity $|v-v_1|$, i.e., $q(v-v_1,\sigma) = q(\sigma)$, but the algorithm below can be  modified to apply to a general set of collision kernels~\cite{bird1970direct,pareschi2001introduction}. Under this assumption,~\eqref{eq:homoBoltz} can be rewritten in the form of 
\begin{equation}\label{eq:homoBoltz_DSMC_1}
    \frac{\partial f}{\partial t} =  [P(f,f) - \mu f],
\end{equation}
where $\rho = \int_\Rthree f dv$, $\mu = \rho \int_\Stwo q(\sigma)d\sigma$, and
\begin{equation*} 
    P(f,f) = \int_\Rthree\int_\Stwo q(\sigma) f'f'_1  d\sigma dv_1.
\end{equation*}
We remark that $f^\rho = \frac{1}{\rho}f$ is a probability density in the velocity space for any given $t$. It also follows the Boltzmann equation~\eqref{eq:homoBoltz} with a scaled collision kernel. Without loss of generality, we regard $f = f^\rho$ as a probability density function hereafter.

In the DSMC method, we consider a set of $N$ velocities evolving in time due to collisions whose distribution can be described by the probability distribution function $f$ in~\eqref{eq:homoBoltz}. We divide time interval $[0,T]$ into $M$ number of sub-intervals of size $\Delta t$. At the $k$-th time interval, the particle velocities are represented as
\begin{equation}~\label{eq:vk}
V_k =\{v_1, \ldots , v_N\} (t_k)
\end{equation}
and we denote the $i$-th velocity particle in $V_k$ as $v_i(t_k)$ or $v_{k,i}$.
The forward Euler scheme applied to~\eqref{eq:homoBoltz_DSMC_1} gives
\begin{equation}\label{eq:homoBoltz_DSMC_2}
    f_{k+1} = (1-\mu\Delta t) f _k + \mu\Delta t\dfrac{P(f _k,f _k)}{\mu},
\end{equation}
where $f _k = f(v,k\Delta t)$. 
If additionally we discretize the velocity distribution at time $t_k = k\Delta t$ by the velocity particles $V_k$ defined in~\eqref{eq:vk}, which is to say
\begin{equation*} 
f _k = f(v,t_k) = f(v,k\Delta t) \approx \frac{1}{N}\sum_{i=1}^N \delta(v-v_{k,i}),\quad k=1,2,\dots,M,
\end{equation*}
we can interpret~\eqref{eq:homoBoltz_DSMC_2} in terms of probability, which is the core idea of the method. 

At time $t_k$, a particle with velocity $v_{k,i}$ will not collide with probability $(1-\mu\Delta t)$, and it will collide with another velocity particle with probability $\mu\Delta t$, according to the collision law described by $P(f _k,f _k)(v)$. Nanbu proposed an algorithm based on this probabilistic interpretation~\cite{nanbu1980direct},  and later its convergence was proved by Babovsky and Illner~\cite{babovsky1989convergence}. One can view~\Cref{alg:DSMC} as the realization of DSMC with the forward Euler scheme~\eqref{eq:homoBoltz_DSMC_2}.

\begin{algorithm}[ht!]
  \caption{Nanbu--Babovsky Algorithm for Maxwellian Molecules\label{alg:DSMC}}
\begin{algorithmic}[1]
\State Compute the initial velocity of particles based on the given initial condition, $V_0 = \{v_{0,1},\dots,v_{0,N}\}$. Set $N_c =\ceil[\big]{N\Delta t \mu}$ and $M=T/\Delta t$ for final time $T$.
\For{$k=1$ to $M$}
\State Given the velocity of particles from the previous time step, $V_{k-1}$.
\State Select $N_c/2$  collision pairs $(i,j)$ uniformly among all possible pairs without replacement.
\State For those selected pairs, perform the collision between $v_{k,i}$ and $v_{k,j}$ based on~\eqref{BoltzmannSolution}; obtain the post-collision velocity ${v}'_{k,i}$ and ${v}'_{k,j}$.
\State Set $v_{k+1,i} = {v}'_{k,i}$ and $v_{k+1,j} = {v}'_{k,j}$.
\State Set $v_{k+1,i} = {v}_{k,i}$ for all particles that have not collided.
\EndFor
\end{algorithmic}
\end{algorithm}

\begin{remark}
The decomposition~\eqref{eq:homoBoltz_DSMC_1} is convenient to illustrate the core idea of DSMC but does not apply to all collision models. The acceptance-rejection method using ``virtual collisions'' is required to sample from a more general cross section ~\cite{caflisch1998monte}.
\end{remark}



\section{Continuous adjoint Boltzmann equations}\label{ContinuousFormulation}
We consider an idealized optimization problem for the spatially homogeneous Boltzmann equation~\eqref{eq:homoBoltz}. 
The initial condition is
\begin{equation}
f(v,0)=f_0(v;\alpha), \label{BoltzmannIC}
\end{equation}
in which $f_0$ is the prescribed initial data depending on a parameter $\alpha$. We aim to find $\alpha$ which optimizes the objective function at time $t=T$,
\begin{equation}\label{eq:OTD_obj}
J_1 (\alpha) = \int_{\bbR^3} r(v) f(v,T) dv.
\end{equation}

When the number of unknowns is large, which is the dimensionality of $\alpha$ in our case, the adjoint-state method is necessary as an efficient numerical method for computing the gradient of a function or operator in a numerical optimization problem~\cite{biegler2003large,cao2003adjoint}. It has applications in geophysics~\cite{chavent1975history}, seismic imaging~\cite{plessix2006review} and general inverse problems~\cite{mcgillivray1990methods}. It is also the theoretical foundation that gave rise to the back-propagation method in the 1990s for neural networks and machine learning~\cite{lecun1988theoretical}.
One outstanding advantage of the adjoint-state method is that the number of PDE solves is independent of the dimension of the parameter for which the gradient needs to be calculated. 

\subsection{Derivation of  the continuous adjoint equation}  

Following the adjoint-state method, we aim to derive the adjoint equations for our optimization problem starting with the Lagrangian
\begin{equation*}
J=\underbracket{\int_\Rthree r(v) f(v,T) dv}_{J_1} +\underbracket{\int_\Rthree \kappa(v) (f(v,0)-f_0(v;\alpha)) dv}_{J_2} +\underbracket{\int_0^T\int_\Rthree \gamma(v,t) ( \partial_t  f(v,t) - Q(f,f) ) dv dt}_{J_3}
\end{equation*}
in which $\kappa(v)$ in $J_2$ is a Lagrange multiplier that enforces the initial condition for any $v\in \bbR^3$, and $\gamma(v,t)$  in $J_3$ is a Lagrange multiplier that enforces the Boltzmann equation for any  $v$ and $t$. 

We remark that $f(v,t)$ is any function here, and its dependence on the Boltzmann equation and the initial condition is imposed through Lagrange multipliers.

We rewrite $J_3$ to calculate its Fr\'echet derivative with respect to $f$.
\begin{equation*}
J_3 = \underbracket{\int_0^T \int_{\bbR^3} \gamma(v,t) \partial_t  f(v,t)  dv dt}_{J_{31}} + \underbracket{(-1) \int_0^T \int_{\bbR^3} \gamma(v,t)  Q(f,f) ) dv dt}_{J_{32}}.
\end{equation*}
After integrating by parts, we have
\begin{equation*}
J_{31}= - \int_0^T \int_{\bbR^3}  f(v,t)  \partial_t  \gamma(v,t)  dv dt + \int_{\bbR^3} \gamma(v,T)   f(v,T) dv - \int_{\bbR^3} \gamma(v,0 )  f(v,0) dv,
\end{equation*}
and for $0<t<T$ the Fr\'echet derivatives are 
\begin{equation*}
\frac{\delta J_{31} }{ \delta {f(v,t)}} = - \partial_t  \gamma(v,t),\quad \frac{\delta J_{31} }{\delta {f(v,T)} } = \gamma(v,T),\quad 
 \frac{\delta J_{31} }{ \delta {f(v,0)} } = -\gamma(v,0 ).
\end{equation*}

The calculation of $ \frac{\delta J_{32} }{ \delta f} $ follows the derivation of collisional invariants for the Boltzmann equation\cite{cercignani1988boltzmann}. 
It follows, using the symmetries~\eqref{eq:symmetry}, that
\begin{align*}
J_{32}&=- \int_0^T \int_\Rthree \gamma(v,t)  Q(f,f)   dv dt = -\int_0^T \int_\Rthree \int_\Rthree \int_\Stwo  \gamma  (f_1' f' - f_1 f) q  d\sigma dv_1 dv dt  \nonumber \\
&=- \int_0^T \int_\Rthree \int_\Rthree \int_\Stwo  (  \gamma' - \gamma )   f_1 f q  d\sigma dv_1 dv dt \quad \mbox{(by switching $v$, $v'$ in $\gamma  f_1' f' $)} \nonumber  \\
&=-  \int_0^T \int_\Rthree \int_\Rthree \int_\Stwo  (  \gamma_1' - \gamma_1 )   f_1 f q  d\sigma dv_1 dv dt \quad \mbox{(by switching $v$, $v_1$)}  \nonumber \\
&=-\frac{1}{2}  \int_0^T  \int_\Rthree \int_\Rthree \int_\Stwo (  \gamma_1' + \gamma' - \gamma_1 - \gamma)   f_1 f q  d\sigma dv_1 dv dt, 
\end{align*}
where $\gamma,\gamma_1,\gamma'$ and $\gamma_1'$ are shorthand notations for $\gamma(v,t),\gamma(v_1,t),\gamma(v',t)$ and $\gamma(v_1',t)$.

Now we perturb $f$ by an amount $\delta f$ and investigate the resulting change in $J_{32}$. We obtain the first-order variation:
\begin{align*}
\delta J_{32}&=-\frac{1}{2} \int_0^T \int_\Rthree \int_\Rthree \int_\Stwo(  \gamma_1' + \gamma' - \gamma_1 - \gamma)  ( f_1 \delta f+ f \delta f_1 ) q  d\sigma dv_1 dv dt   \nonumber \\
&=-    \int_0^T\int_\Rthree \int_\Rthree \int_\Stwo(  \gamma_1' + \gamma' - \gamma_1 - \gamma)   f_1 \delta f  q  d\sigma dv_1  dv dt. 
\end{align*}
Thus, the Fr\'echet derivative of $J_{32}$ with respect to $f$ is
\begin{equation*}
\frac{\delta J_{32} }{ \delta f }= -    \int_\Rthree \int_\Stwo  (  \gamma_1' + \gamma' - \gamma_1 - \gamma)   f_1   q  d\sigma dv_1  =  -L^*[f] \gamma,\end{equation*}
where the second equality follows~\eqref{eq:adj_linear_operator}.
Together with all the other functional derivatives, we obtain the following equations:
\begin{align*}
\frac{\delta J}{ \delta {f(v,t)}} &= - \partial_t  \gamma  - L^*[f] \gamma,\\
\frac{\delta J}{ \delta {f(v,T)}} &=\gamma(v,T) + r(v), \\
\frac{\delta J}{ \delta {f(v,0)} } &=  - \gamma(v,0) +\kappa(v),  \\
\frac{\partial J}{ \partial{\alpha} } &= -\int_\Rthree \kappa(v) \partial_\alpha f_0(v;\alpha) dv. 
\end{align*}

Based on the first-order necessary condition for optimality, we set these variations to zero to obtain the adjoint equation for the adjoint variable $\gamma$:
\begin{subequations}~\label{eq:backwards}
\begin{align}
- \partial_t  \gamma &=  L^*[f] \gamma, \label{backwards1} \\
\gamma(v,T) &= - r(v).  \label{backwards0}
\end{align}
\end{subequations}
The adjoint equation evolves backward in time from $T$ to $0$, thus the final condition~\eqref{backwards0} is given instead of the initial condition. After time evolution, we can use $\gamma(v,0)$ to eliminate $\kappa(v)$ and obtain the gradient of the objective function with respect to parameter $\alpha$:
\begin{equation}  \label{optimization1}
\partial_\alpha J = -\int_\Rthree \gamma(v,0) \partial_\alpha f_0(v;\alpha) dv.
\end{equation}
The continuous adjoint equation~\eqref{backwards1} is an integro-differential equation which shares significant similarities with the linearized Boltzmann Equation~\cite{cercignani1969mathematical}. We will further analyze and discuss the adjoint system in~\Cref{sec:discussion}.

Given an initial guess $\alpha_0$ for the model parameter, the first step is to solve the forward equations (\ref{eq:homoBoltz}) numerically with initial condition (\ref{BoltzmannIC}) for $f$. A common choice is the DSMC method because of the nonlinear collision operator. The second step is to solve the adjoint equation~\eqref{backwards1} numerically with the final condition~\eqref{backwards0} for the adjoint variable $\gamma$. One can solve the adjoint equation using a finite difference method in $t$ and grid-based quadrature in $v$, but this is computationally difficult because $v$ is three dimensional, and the integral at each value of $v$ is five-dimensional. In~\Cref{sec:DSMCtype_method}, we propose~\Cref{alg:adj_DSMC} as an efficient Monte Carlo type method for solving~\eqref{backwards1}. In the third step, we compute the gradient $\partial_{\alpha} J$ based on~\eqref{optimization1} using the numerical solutions from the first two steps. Numerical examples and comparisons are presented in~\Cref{sec:numerics}. One can then update $\alpha$ iteratively using a gradient-based optimization algorithm to find the optimal parameter. This is the ``optimize-then-discretize'' (OTD) approach to solve the PDE-constrained optimization numerically. 

We remark that the entire derivation above is unrelated to the dimensionality of the model parameter $\alpha$. In practice, $\alpha$ can be, for example, a function of the velocity domain.  Hence, an accurate numerical discretization of the parameter can contain thousands of variables. It is then infeasible to obtain the gradient by numerical differentiation. On the other hand, using the adjoint-state method, the cost of evaluating the gradient once is equivalent to only one forward solve of the Boltzmann equation~\eqref{eq:homoBoltz} and one (backward) solve of the adjoint equation~\eqref{backwards1}, independent of the number of variables in $\alpha$. The adjoint-state method is an extremely efficient tool for large-scale optimization problems.



\subsection{A particle method for the continuous adjoint Boltzmann equation}  \label{sec:DSMCtype_method}
The adjoint DSMC method presented in~\Cref{DiscreteFormulation} is the main result of this work and the most efficient and accurate that we have found for solving the adjoint problem. Nevertheless, the particle method described in this section may be of independent interest.  It is based on a formal derivation with several steps that are not fully justified, but are validated by the numerical results  in~\Cref{sec:numerics}.

The continuous adjoint equation~\eqref{backwards1} is based on the solution to the forward equation $f(v,t), t\in [0,T]$. Therefore, we assume that the state (forward) equation~\eqref{eq:homoBoltz} has been solved numerically by the DSMC method, where the choice of collision times, collision partners, and collision angles during the simulation are all stored in memory. Based on the solution from the forward DSMC simulation, we will (approximately) solve the equation~\eqref{backwards1} for $\gamma(v,t)$ backward in time $t$ with the ``final data"~\eqref{backwards0}.

According to~\eqref{eq:homoBoltz} and~\eqref{backwards1}, we have
\begin{equation*}
\partial_t (\gamma f) = \gamma \partial_t f + f \partial_t \gamma = \iint(f'f_1' - f f_1) \gamma qd\sigma dv_1 -  \ \iint (\gamma_1' + \gamma' - \gamma_1 -\gamma) f_1 f q  d\sigma dv_1. 
\end{equation*}
We then multiply both sides by a time-independent test function $\psi(v)$ and integrate over the velocity domain. The left-hand side becomes
\begin{equation}\label{eq:LHS_test}
\int_\Rthree \psi \partial_t (\gamma f) dv = \partial_t \left(  \int_\Rthree \psi \gamma f dv\right).
\end{equation}
Using the symmetries~\eqref{eq:symmetry}, the right-hand side becomes
\begin{align}
&\quad \iiint \psi \bigg\{(f'f_1' - f f_1) \gamma - (\gamma_1' + \gamma' - \gamma_1 -\gamma) f_1 f \bigg\} qd\sigma dv_1 dv \nonumber \\
&= \frac{1}{2} \iiint
\left(\psi_1'\gamma_1' + \psi'\gamma' -\psi_1\gamma_1 - \psi\gamma \right) f_1 f qd\sigma dv_1 dv  - \frac{1}{2} \iiint \left(\psi + \psi_1  \right)(\gamma_1' + \gamma' - \gamma_1 -\gamma) f_1 f qd\sigma dv_1 dv \nonumber \\
&= \frac{1}{2} \iiint \bigg\{ 
\psi_1'\gamma_1' + \psi'\gamma' - \psi_1(\gamma_1' + \gamma' -\gamma)
-\psi (\gamma_1' + \gamma' -\gamma_1)
\bigg\} f_1 f qd\sigma dv_1 dv \label{eq:RHS_test}.
\end{align}

Combining both~\eqref{eq:LHS_test} and~\eqref{eq:RHS_test}, we discretize the time domain of the equation using the Euler scheme in $t$. Consider a time interval $[t_k, t_{k+1}]$ during which there is a collision of $v$ and $v_1$ at $t=t_k$ to produce $v'$ and $v_1'$ at $t=t_{k+1}$.  

As in the (forward) DSMC~\Cref{alg:DSMC},   $f=f(v,t)$ and $f_1=f(v_1,t)$ in~\eqref{eq:RHS_test} should be evaluated at $t_k$ since $v$ and $v_1$ are independent before the collision so that the product $f f_1$ is the joint density function for $v$ and $v_1$. 

On the other hand, the choice of time $t$ at which to evaluate $\gamma(v,t)$ is a new issue and not so clearcut.  We choose to evaluate $\gamma(v,t)$  at $t_{k+1}$ based on the following considerations: 
First, the data for the adjoint variable comes at the final time~\eqref{backwards0}, so that~\eqref{backwards1} should be solved backward in time, and it is most natural for the right-hand side~\eqref{eq:RHS_test} to be evaluated at the final time $t_{k+1}$ of the time interval.
Second,  since $\gamma$ is solved backward in time, we expect that 
 $\gamma(v,t_{k+1})$ is independent of $f(v,t_k)$, so that it is reasonable to sample using these values. 
Third, our numerical computations in~\Cref{sec:numerics} verifies that this choice leads to a correct result.


Based on these choices, we take  $\gamma(v,t)=\gamma_{k+1}(v)$ at $t=t_{k+1}$ and $f(v,t)=f_k(v)$ at $t=t_k$ in~\eqref{eq:RHS_test} to obtain (with $\psi$ independent of $t$)
\begin{align}
\int_\Rthree \psi \gamma_{k+1} f_{k+1} dv -\int_\Rthree \psi \gamma _k f _k dv 
&\approx\frac{\Delta t}{2} \iiint (\psi_1'\gamma_{k+1,1}' + \psi'\gamma_{k+1}')  f_k(v_1) f_k(v) q(\sigma)d\sigma dv_1 dv  \nonumber \\
 &\quad-\frac{\Delta t}{2}  \iiint  \psi_1(\gamma_{k+1,1}' + \gamma_{k+1}' -\gamma_{k+1}) f_k(v_1) f_k(v) q(\sigma)d\sigma dv_1 dv  \nonumber\\
 &\quad -\frac{\Delta t}{2}  \iiint\psi (\gamma_{k+1,1}' + \gamma_{k+1}' -\gamma_{k+1,1})f_k(v_1) f_k(v) q(\sigma)d\sigma dv_1 dv. \label{eq:OTD_DSMC}
\end{align}

Without loss of generality, we assume $\rho = \int  f dv = 1$ and 
$\mu = \rho \int q(\sigma) d\sigma = 1$. Then $f(v)$ is a probability density in $\Rthree$, and $F(\sigma,v,v_1) =  f_k(v) f_k(v_1)q(\sigma)$ is a probability density function in the product space $\Stwo \times \Rthree \times \Rthree$. We apply Monte Carlo quadrature to approximate the integrals on both sides of~\eqref{eq:OTD_DSMC}, using the  velocities from  $V_{k}$ and $V_{k+1}$ (defined as in~\eqref{eq:vk}) and values of $\sigma$ that were selected in the forward DSMC calculation. Similar to~\eqref{eq:vk}, we represent the adjoint variables as
\begin{equation*} 
\Gamma_{k+1} =\{\hat \gamma_1, \ldots , \hat \gamma_i,\ldots \hat \gamma_N\} (t_{k+1}).
\end{equation*}

Both sides of the resulting Monte Carlo sums for~\eqref{eq:OTD_DSMC} are  nonzero only for velocities that undergo collisions. For the  Monte Carlo quadrature we have  $N_c/2 =  \ceil[]{ \mu\Delta t  N}/2$ collision pairs. 

The $j$-th collision between velocities $v_{j},v_{j_1} \in V_k$, with collision parameters $\sigma_{j}$,   results in velocities $v_{j}',v_{j_1}' \in V_{k+1}$. The corresponding values of $\gamma$ are denoted as $\hat \gamma, \hat \gamma_1 \in \Gamma_k$  and $\hat \gamma', \hat \gamma_1' \in \Gamma_{k+1}$. 

On the left-hand side of~\eqref{eq:OTD_DSMC}, the term $\gamma_{k+1}$ is represented by $\hat \gamma', \hat \gamma_1'$ in $\Gamma_{k+1}$ and the term $\gamma_{k}$ is represented by $\hat \gamma, \hat \gamma_1 \in \Gamma_k$. 

On the right-hand side of~\eqref{eq:OTD_DSMC}, the terms $\gamma_{k+1}',\gamma_{k+1,1}'$ are represented by $\hat \gamma', \hat \gamma_1' \in\Gamma_{k+1}$, but  the terms $\gamma_{k+1},\gamma_{k+1,1}$ are at $t=t_{k+1}$ and correspond to velocities $v, v_1 \notin V_{k+1}$, so that the corresponding $\hat \gamma$ values are not in  $ \Gamma_{k+1}$. We will instead use the values of the continuous function $\gamma(v,t_{k+1})$ and $\gamma(v_1,t_{k+1})$.

The resulting quadrature approximation for~\eqref{eq:OTD_DSMC} by Monte Carlo sampling is
\begin{align} \label{eq:OTD_DSMC2}
&\quad \sum_{j=1}^{ N_c/2} \psi'_{j_1} \hat \gamma_{j_1}(t_{k+1}) + \psi'_{j} \hat \gamma_{j}(t_{k+1}) - \psi_{j_1} \hat \gamma_{j_1}(t_{k})  -  \psi_{j} \hat \gamma_{j}(t_{k}) \\ 
&\approx \sum_{j=1}^{ N_c/2}  \psi'_{j_1} \hat \gamma_{j_1}(t_{k+1}) +  \psi'_{j}  \hat \gamma_{j}(t_{k+1})  -  \psi_{j_1} \left( \hat \gamma_{j_1}(t_{k+1})+ \hat \gamma_{j}(t_{k+1}) - \gamma(v_j, t_{k+1})\right) \nonumber\\
& - \psi_j \left( \hat \gamma_{j_1}(t_{k+1}) + \hat \gamma_{j}(t_{k+1}) - \gamma(v_{j_1} , t_{k+1}) \right), \nonumber
\end{align}
where $\psi_j = \psi(v_j)$, $\psi_{j_1} = \psi(v_{j_1})$, $\psi'_j = \psi(v'_j)$ and $\psi'_{j_1} = \psi(v'_{j_1})$.

The values of  the continuous adjoint function in~\eqref{eq:OTD_DSMC2} can be expressed in terms of a conditional expectation: 
\begin{equation*}
    \gamma(v,t)= \mathbb E [\hat \gamma_i(t)|v_i = v].
\end{equation*}

Since $\psi$ is an arbitrary test function, and~\eqref{eq:OTD_DSMC2} holds for any $\psi$, we can match the coefficients and obtain a Monte Carlo type numerical scheme for solving $\gamma(v,t)$.
If velocity particle $v_j$ did not collide at time $t_k$, 
\begin{equation}
   \hat \gamma_j(t_k) = \hat \gamma_j(t_{k+1}).
\end{equation}
If $v_j$ and $v_{j_1}$ are a collision pair at time $t_{k}$,
\begin{subequations}\label{eq:OTD_DSMC3}
\begin{align}
\hat \gamma_j(t_k) &= \hat \gamma_j(t_{k+1}) + \hat \gamma_{j_1}(t_{k+1}) - \mathbb E[\hat \gamma_i (t_{k+1}) |v_i = v_{j_1}],\\
\hat \gamma_{j_1}(t_k) &= \hat \gamma_{j}(t_{k+1}) + \hat \gamma_{j_1}(t_{k+1}) -\mathbb E[\hat \gamma_i(t_{k+1}) |v_i = v_j].
\end{align}
\end{subequations}

Following~\cite{longstaff2001valuing,TsitsiklisRoy},
one may approximate the expectation $\mathbb E[\hat \gamma_i (t_{k+1}) |v_i = v_{j}]$ by numerical interpolation as the following:
\begin{equation}\label{eq:OTD_DSMC4}
\hat \gamma(v_{j} , t_{k+1}) = \mathbb E[\hat \gamma_i (t_{k+1}) |v_i= v_{j}] \approx \frac{1}{L}\sum_{l=1}^L \omega(v_{j}-v_l) \hat \gamma_l (t_{k+1})
\end{equation}
for appropriately chosen interpolation coefficients $\omega(v)$.

Given the final condition~\eqref{backwards0}, the continuous adjoint equation~\eqref{backwards1} can be solved numerically following Algorithm~\ref{alg:adj_DSMC}.
\begin{algorithm}[ht!]
  \caption{Algorithm for solving the continuous adjoint equation~\eqref{eq:backwards}\label{alg:adj_DSMC}}
\begin{algorithmic}[1]
\State Given the final-time velocity particles $V_{M}$ in the forward DSMC and the final condition~\eqref{backwards0}, set $\hat \gamma_i(T) = -r(v_{M,i}),\ i =1,\dots,N.$ Obtain $\Gamma_M$.
\For{$k=M-1$ to $0$}
\State Given $\Gamma_{k+1}$ from the previous iteration and $V_{k}$ from the forward DSMC.
\If {$v_j\in V_{k}$ did not collide at $t_k$}
\State Set $\hat \gamma_j(t_{k}) = \hat \gamma_j(t_{k+1})$.
\ElsIf {$v_j, v_{j_1}\in V_{k}$ collided at $t_k$}
\State Approximate $\mathbb E[\hat \gamma_i (t_{k+1}) |v_i = v_j]$ and $\mathbb E[\hat \gamma_i (t_{k+1}) |v_i= v_{j_1}]$ following~\eqref{eq:OTD_DSMC4}.
\State Set $\hat \gamma_j(t_k)$ and $\hat \gamma_{j_1}(t_k)$ following~\eqref{eq:OTD_DSMC3}. 
\EndIf
\State Obtain $\Gamma_k$.
\EndFor
\end{algorithmic}
\end{algorithm}
\begin{remark}
We regard the numerical method in~\Cref{alg:adj_DSMC} as the DSMC-type scheme for the continuous adjoint equation,
depending on a solution of~\eqref{eq:homoBoltz}  by the forward DSMC method following~\Cref{alg:DSMC}. 
The gradient of the objective function~\eqref{optimization1} can be computed using the values  of $\gamma(v,t)$ at $t=0$, i.e., $\Gamma_0$. Numerical examples are shown in~\Cref{sec:numeric_DSMCtype_method}.
\end{remark}




\section{Adjoint DSMC for Boltzmann equation} \label{DiscreteFormulation}
In this section, we derive the adjoint system based on the DTO approach. We first rewrite the objective functions and the constraints based on the DSMC method, to which the Lagrange multiplier method will be applied.

\begin{figure} 
\begin{center}
\begin{tikzpicture}
\centering
\begin{scope}
\draw[black, thick] (0,10) -- (5,10);
\draw[black, thick] (1,10) -- (2,9);
\draw[black, dashed, thick] (2,9) -- (3.3,8);
\draw[black, dashed, thick] (3,9.5) -- (2,9);
\draw[black, thick] (2,9) -- (1.5,8);
\node (a) at (1,9) {$w_F,\bgamma_F$};
\node (a) at (5.5,10) {$t = T$};
\node (a) at (2,7) {$(a)$};
\end{scope}

\begin{scope}[shift={(-1,0)}]
\draw[black, thick] (8,8) -- (13,8);
\draw[black, thick] (9,8) -- (12,11);
\draw[black, thick] (10.5,9.5) -- (9,11);
\draw[black, thick] (10.5,9.5) -- (12,8.5);
\node (a) at (9,10) {$w_{k+1},\bgamma_{k+1}$};
\node (a) at (12,10) {$\tilde{w}_{k+1},\tilde{\bgamma}_{k+1}$};
\node (a) at (9,9) {$w_I,\bgamma_I$};
\node (a) at (12,9) {$\tilde{w}_{k},\tilde{\bgamma}_{k}$};
\node (a) at (13.5,8) {$t = 0$};
\node (a) at (11,7) {$(b)$};
\end{scope}
\end{tikzpicture}
\caption{$(a)$~Collision for velocity particles at the final time $t=T$; $(b)$~Collision for velocity particles at the initial time $t=0$.}~\label{fig:collision1}
\end{center}
\end{figure}
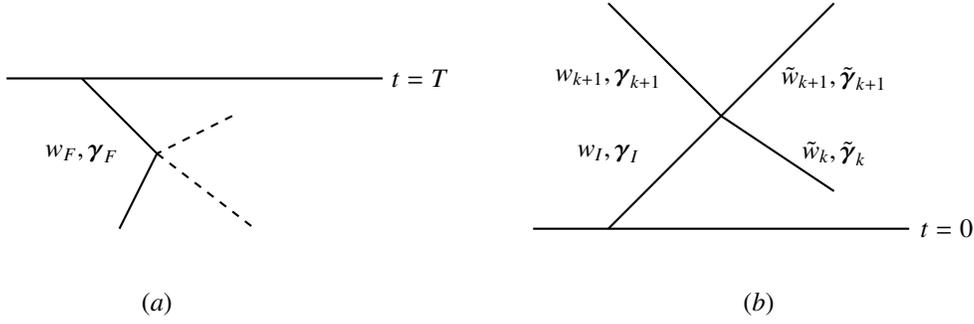
Following Algorithm~\ref{alg:DSMC}, we discretize the time $[0,T]$ into $M$ equal intervals $[t_k,t_{k+1}]$ and denote $v_{k,i}$ to be the velocity of the $i$-th particle at time $t_k$, for $1\leq i\leq N$ and $0\leq k\leq M$. We also denote $v_{I,i}$ to be the velocity of the $i$-th particle at the initial time $t=0$, and $v_{F,i}$ to be the velocity of the $i$-th particle at the final time $t=T$. The DSMC method applied to the spatially homogeneous Boltzmann equation consists of a sequence of particle collisions in sequential order, and the time step $\Delta t$ plays no role. So we will choose $\Delta t$ in whatever way is most convenient for the exposition. In this section, $\Delta t$ will be chosen small enough that there is at most a single collision in each time step.

Consider the new objective function $\cJ$ defined as
\begin{equation} \label{eq:DTO_obj}
\cJ = \underbracket{\frac{1}{N}\sum_{i=1}^N r(v_{{F},i})}_{\cJ_1} +\underbracket{ \frac{1}{N}\sum_{i=1}^N \bgamma_{I,i} \cdot (v_{I, i}- v_{0,i}(\alpha) )}_{\cJ_2}+\underbracket{\frac{1}{N} \sum_{k=1}^M \sum_{i=1}^N \bgamma_{k,i} \cdot (v_{k+1,i}-v_{k,i}')}_{\cJ_3}.
\end{equation}
Here, $\cJ_1$ is the Monte Carlo quadrature of the objective function~\eqref{eq:OTD_obj} by particle velocities $\{v_{F,i}\}_{i=1}^N$, $\cJ_2$ is the constraint on the DSMC initial condition $v=v_0(\alpha)$ using the Lagrange multiplier $ \{\bgamma_{I,i}\}_{i=1}^N$, and $\cJ_3$ is the constraint that enforces the binary collision law~\eqref{BoltzmannSolution} using the Lagrange multiplier $\bgamma_{k,i}$ for each particle $i$ at the $k$-th time interval. In particular, $v_{k,i}'$ represents the post-collision velocity of particle $i$ if it participates in the collision at the $k$-th time interval. Otherwise, $v_{k,i}' = v_{k,i}$, which means the particle velocity remains the same at the $(k+1)$-th time interval.

Again, we regard that \{$v_{k,i}$\} is a general set of velocities, and its dependence on the collision rules~\eqref{BoltzmannSolution} is imposed through the Lagrange multipliers.
The following variation is easily calculated:
\begin{equation}
\partial_\alpha \cJ=- \frac{1}{N} \sum_{i=1}^N \bgamma_{{I},i} \cdot \partial_\alpha v_{0,i}(\alpha) \label{DalphaEqtn}
\end{equation}

We proceed to derive the equations in the backward time direction from the final time $t=T$ to the initial time $t=0$. First, we consider the equations that come from the derivative of $\cJ$ with respect to $v_{{F},i}$, which are derived in \Cref{DiscreteFinal}. Second, we consider collisions that occur at a time $t$ where $0 < t<T$. 
The derivatives of $\cJ$ with respect to the velocities in these collisions are derived in \Cref{DiscreteInterior}.

\begin{figure}
\begin{center}
\begin{tikzpicture}[scale=0.95]
\centering
\draw[black, thick] (1,.2) -- (3,2);
\draw[black, thick] (3,2) -- (7,4);
\draw[black, thick] (3,2) -- (5,.2);
\draw[black, dashed, thick] (2,3) -- (3,2);
\draw[black, thick] (7,4) -- (10,7);
\draw[black, thick] (7,4) -- (4,7);
\draw[black, thick] (7,4) -- (10,2);
\draw[black, thick] (8,.2) -- (10,2);
\draw[black, thick] (10,2) -- (12,.2);
\draw[black, dashed, thick] (10,2) -- (11.5,2.8);
\node (a) at (0.8,1) {$\tilde w_{k-1},\tilde \bgamma_{k-1}$};
\node (a) at (5.2,1) {$w_{k-1}, \bgamma_{k-1}$};
\node (a) at (4,3.1) {$w_{k},\bgamma_{k}$};
\node (a) at (4,6) {$w_{k+1},\bgamma_{k+1}$};

\node (a) at (7.8,1) {$\overline{w}_{k-1},\overline{\bgamma}_{k-1}$};
\node (a) at (12.2,1) {$\widetilde{\overline{w}}_{k-1},\widetilde{\overline{\bgamma}}_{k-1}$};
\node (a) at (10.1,3.1) {$\tilde{w}_{k} ,\tilde{\bgamma}_{k}$};
\node (a) at (10.2,6){$\tilde{w}_{k+1},\tilde{\bgamma}_{k+1}$};
\end{tikzpicture}
\caption{Collision for velocity particles at intermediate time.}~\label{fig:collision2}
\end{center}
\end{figure}
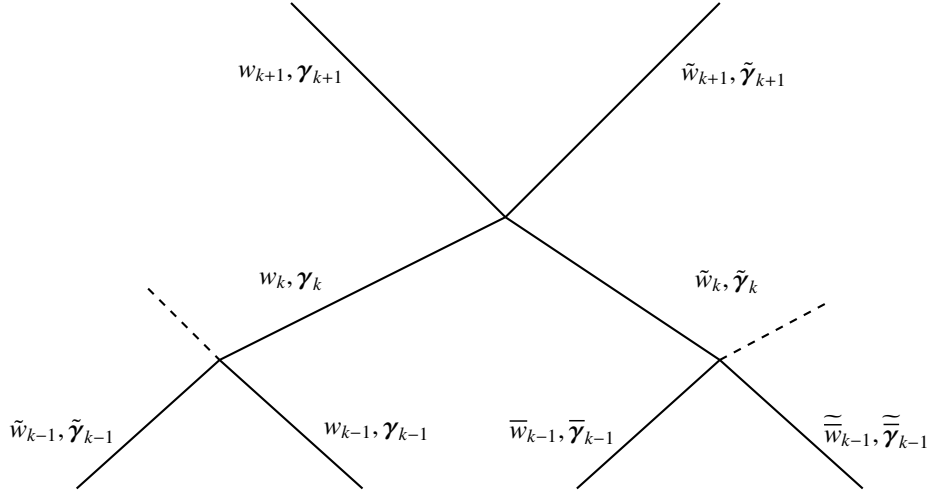


\subsection{Velocities at  the final time for DSMC} \label{DiscreteFinal}
We denote a particular particle velocity  at the final time as $w_F=v_{{F}, i}$ for some $i$, and the corresponding value of its dual variable is  $\bgamma_F = \bgamma_{{F}, i} $, as shown in Fig.~\ref{fig:collision1}(a). The velocity $w_F$ occurs in the objective function $\cJ$ at two places. First, it directly appears in the term $r(v_{{F}, i})$ in $\cJ_1$. Second, it occurs in the last collision that involves $v_{{F}, i}$, through the term $\bgamma_{F,i} \cdot (v_{F,i}-v_{k,i}')$ in $\cJ_3$. 
In this term,  $v_{k,i}'$ is a velocity produced by the collision, and is a function of the two velocities before the collision, but $v_{k,i}'$ does not depend on $v_{F,i}$ here. Therefore,
\begin{equation*}
\partial_{ w_F } \cJ = \frac{1}{N} \partial_{ w_F } r(w_F ) +\frac{1}{N}  \partial_{ w_F} \left(  \bgamma_{F} \cdot(w_F-v_{k,i}')\right) \nonumber
= \frac{1}{N}\partial_v r(w_F) +\frac{1}{N} \bgamma_{F}.
\end{equation*}
The resulting equation comes from the first-order optimality $\partial_{w_F } \cJ=0$:
\begin{equation}
\bgamma_{F} = - \partial_v r(w_F) \label{gammaFinal}.
\end{equation}
Equation~\eqref{gammaFinal} provides the starting values for solving the adjoint variable $\bgamma$ backward in time.



\subsection{Collisions away from the final time for DSMC} \label{DiscreteInterior}
According to the DSMC method, with $V_k$ denoting the set of velocities in the $k$-th time interval as in (\ref{eq:vk}), the collision operator $C$ takes $V_k$ to $V_{k+1}=C[V_k]$ as follows: We first choose indices $i$ and $j$, and denote the corresponding  velocities as $(w_k,{\tilde w}_k) = (v_{k,i}, v_{k,j})$ from $V_k$ and also choose a unit vector $\sigma_k$. Here, we assume the $i$-th and the $j$-th particles collide at the $k$-th time interval. Then $C$ takes $(w_k,{\tilde w}_k)$ to
\begin{equation*}
(w_{k+1},{\tilde w}_{k+1}) = (w_k',{\tilde w}_k')
\end{equation*}
where
\begin{subequations}\label{eq:w_eqn_both}
\begin{align}
w_k' &=\frac{1}{2} ( (w_k+{\tilde w}_k) + |w_k- {\tilde w}_k | \sigma_k )\label{eq:w_eqn}, \\
{\tilde w}_k' &=\frac{1}{2} ( (w_k+{\tilde w}_k) - |w_k- {\tilde w}_k | \sigma_k ). \label{eq:w_tilde_eqn} 
\end{align}
\end{subequations}
Note that the resulting velocities $(w_{k+1},{\tilde w}_{k+1})$ correspond to $(v_{k,i}', v_{k,j}')$.

On the other hand, we denote $w_{k-1}$ and $\tilde w_{k-1}$ to be the particle velocities whose collision produces $w_k$. Similarly, ${\tilde w}_k$ is produced by particle velocities $\overline w_{k-1}$ and $\widetilde{\overline{w}}_{k-1}$. Particle velocities 
$$w_{k-1},\ \tilde w_{k-1},\ \overline{w}_{k-1},\ \widetilde{\overline{w}}_{k-1} ,\  w_k,\ \tilde w_k, w_{k+1},\ \tilde w_{k+1}$$ 
correspond to the dual (adjoint) variables 
$$
\bgamma_{k-1},\ \tilde\bgamma_{k-1},\ \overline{\bgamma}_{k-1},\  \widetilde{\overline{\bgamma}}_{k-1},\ \bgamma_k,\ \tilde \bgamma_k,\  \bgamma_{k+1},\ \tilde\bgamma_{k+1},
$$ 
respectively, as shown in Fig.~\ref{fig:collision2}, which schematically depicts the three collisions that involve particles $w_k$ and ${\tilde w}_k$.

Next, we recall the variational quantity $\cJ_3$ for collisions, with dual quantities $\bgamma$, as
\begin{equation*}
\cJ_3 = \frac{1}{N}  \sum_{k} \sum_{i} \bgamma_{k,i} \cdot (v_{k+1,i}-v_{k,i}') .
\end{equation*}
We use $\cJ_{3k}$ to denote the terms in $\cJ_3$ that involve $w_k$ or ${\tilde w}_k$: 
\begin{equation*}
\cJ_{3k} = \frac{1}{N} \big\{\bgamma_k \cdot (w_k-w_{k-1}') + {\tilde \bgamma}_k \cdot ({\tilde w}_k-\overline{w}_{k-1}')  \nonumber + \bgamma_{k+1} \cdot (w_{k+1}-w_k')     + {\tilde \bgamma}_{k+1} \cdot ({\tilde w}_{k+1}-{\tilde w}_k')\big\}.
\end{equation*}
Note that $w_{k-1}'$  is a  function only of $w_{k-1}$ and $\tilde w_{k-1}$, and that $\overline{w}_{k-1}'$ is a function only of   $\overline{w}_{k-1}$ and $\widetilde{\overline{w}}_{k-1}$. Consequently, $w_{k-1}'$ and $\overline{w}_{k-1}'$, as well as   $w_{k+1}$ and $\tilde w_{k+1}$, do not depend on $w_k$ and ${\tilde w}_k$. 
It follows that 
\begin{subequations}\label{eq:fwd_perturb45}
\begin{align}
\partial_{w_{k}} \cJ_3 &= \partial_{w_k} \cJ_{3k} = \frac{1}{N}  \left(\bgamma_{k} - \bgamma_{k+1} \cdot \partial_{w_k} w_{k}'  - {\tilde \bgamma}_{k+1} \cdot \partial_{w_k} {\tilde w}_{k}'\right),  \\
\partial_{{\tilde w}_{k}}\cJ_3 &= \partial_{{\tilde w}_k}  \cJ_{3k} = \frac{1}{N} \left({\tilde \bgamma_{k}} - \bgamma_{k+1} \cdot \partial_{{\tilde w}_k} w_{k}'  - {\tilde \bgamma}_{k+1} \cdot \partial_{{\tilde w}_k} {\tilde w}_{k}'\right).
\end{align}
\end{subequations}
These derivatives can be calculated from~\eqref{eq:w_eqn_both} to get
\begin{align*}
\bgamma \cdot \partial_{w_k} w_{k}'  &=  \frac{1}{2}  \left( \bgamma + \bgamma \cdot\sigma_k  {(w_{k} - {\tilde w}_{k} ){\hat ~} }\right)  \nonumber, \\
\bgamma \cdot \partial_{{\tilde w}_{k}} w_{k}'  &=  \frac{1}{2}  \left( \bgamma - \bgamma \cdot\sigma_k  (w_{k} - {\tilde w}_{k} ){\hat ~} \right)  \nonumber, \\
\bgamma \cdot \partial_{w_k} {\tilde w}_{k}'  &= \frac{1}{2}  \left( \bgamma - \bgamma \cdot\sigma_k  { {(w_{k} - {\tilde w}_{k} ){\hat ~}}} \right)  \nonumber, \\
\bgamma \cdot \partial_{{\tilde w}_{k}} {\tilde w}_{k}'  &=  \frac{1}{2}  \left( \bgamma + \bgamma \cdot\sigma_k  { {(w_{k} - {\tilde w}_{k} ){\hat ~}}} \right),
\end{align*}
where $x{\hat ~} = x/|x|$ denotes the unit vector along the direction of $x$. Therefore, 
\begin{align*}
\partial_{w_{k}} \cJ_3 &= \frac{1}{N}  \bigg\{\bgamma_{k} -  \frac{1}{2}( \bgamma_{k+1} + {\tilde \bgamma}_{k+1})  -  \frac{1}{2} ( \bgamma_{k+1} - {\tilde \bgamma}_{k+1}) \cdot\sigma_k  (w_{k} - {\tilde w}_{k} ){\hat ~} \bigg\},  \\
\partial_{{\tilde w}_{k}} \cJ_3 &= \frac{1}{N}  \big\{ {\tilde \bgamma}_{k} -  \frac{1}{2} ( \bgamma_{k+1} + {\tilde \bgamma}_{k+1}) +  \frac{1}{2} ( \bgamma_{k+1} - {\tilde \bgamma}_{k+1}) \cdot\sigma_k  (w_{k} - {\tilde w}_{k} ){\hat ~} \bigg\}.
\end{align*}
Setting these two partial derivatives to $0$, gives the equations
\begin{subequations}\label{gammaEqtns}
\begin{align}
\bgamma_{k} &=    \frac{1}{2}  ( \bgamma_{k+1} + {\tilde \bgamma}_{k+1})  + \frac{1}{2} ( \bgamma_{k+1} - {\tilde \bgamma}_{k+1}) \cdot\sigma_k   (w_{k} - {\tilde w}_{k} ){\hat ~},   \label{gammaEqtn} \\
{\tilde \bgamma}_{k} &=   \frac{1}{2}  ( \bgamma_{k+1} + {\tilde \bgamma}_{k+1}) -  \frac{1}{2} ( \bgamma_{k+1} - {\tilde \bgamma}_{k+1}) \cdot\sigma_k  (w_{k} - {\tilde w}_{k} ){\hat ~}  . \label{gammaTildeEqtn}
\end{align}
\end{subequations}
Moreover,~\eqref{gammaEqtns} can be rewritten using the operator notation~\eqref{eq:AB_def} from Section~\ref{sec:operator_notation} as
\begin{equation}\label{eq:AB_gamma}
        \begin{pmatrix}
    \bgamma_{k+1}\\
    \tilde \bgamma_{k+1}
    \end{pmatrix}     = A(\sigma_k,\alpha_k) \begin{pmatrix}
    \bgamma_{k}\\
    \tilde \bgamma_{k}
    \end{pmatrix}, \quad 
   \begin{pmatrix}
    \bgamma_{k}\\
    \tilde \bgamma_{k}
    \end{pmatrix}     = B(\sigma_k,\alpha_k)    \begin{pmatrix}
    \bgamma_{k+1}\\
    \tilde \bgamma_{k+1}
    \end{pmatrix},
\end{equation}
in which $\alpha_k=(w_{k} - {\tilde w}_{k} ){\hat ~} $.


Note that~\eqref{gammaEqtns} for $\bgamma_{k}$ and ${\tilde \bgamma}_{k}$  does not involve any variables at time interval $k-1$. It follows that~\eqref{gammaEqtns} is valid for values of $\bgamma_I = \bgamma_{k}$ and ${\tilde \bgamma}_{k}$ corresponding to a collision in which one (or both) of the velocities $w_{I}= w_{k}$ started at $t=0$, as depicted in Fig.~\ref{fig:collision1}(b). 

The adjoint DSMC equations~\eqref{gammaEqtn}-\eqref{gammaTildeEqtn}, which we derive following the DTO approach, are the DSMC analog of the dual equation (\ref{backwards1}), and they are one of the main results of this paper. The evolution equations for the dual variable $\bgamma$ are solved backward in time,  from $(k+1)$-th time interval to the $k$-th time interval, for each $k$. For each collision, these equations are solved for the dual variables corresponding to the particle velocity variables involved in that collision.



\subsection{Summary for discrete adjoint system} 

In summary, we solve the direct equations in DSMC forward in time for the velocity particles $V=\{v_1, \ldots , v_N \}$ with the initial particle velocity sampled from~\eqref{BoltzmannIC} given the current value of $\alpha$. Collisions are performed according to~\eqref{BoltzmannSolution}. 
Then the adjoint equations~\eqref{gammaEqtn}-\eqref{gammaTildeEqtn} for the adjoint particle $\bgamma$ are solved backward in time following~\Cref{alg:adj_DSMC1}.
At the initial time $t=0$, we use \eqref{DalphaEqtn} to calculate the gradient of the objective function with respect to the model parameter $\alpha$, which can be used by optimization algorithms to update the $\alpha$ iteratively.

\begin{algorithm}[ht!]
\caption{Algorithm for solving the discrete adjoint DSMC system~\eqref{gammaEqtns}.~\label{alg:adj_DSMC1}}
\begin{algorithmic}[1]
\State Given the final-time velocity particles $V_{M} = \{v_{F,1},\ldots, v_{F,N}\}$ from the forward DSMC, set $\bgamma_{F,j} = - \partial_v r(v_{F,j})$ for $j=1,\ldots,N$ following~\eqref{gammaFinal}.
\For{$k=M-1$ to $0$}
\State Given $\{\bgamma_{k+1,1},\ldots, \bgamma_{k+1,N}\}$ from the previous iteration and collision parameters in the forward DSMC.
\If {$v_j\in V_{k}$ did not collide at $t_k$}
\State Set $\bgamma_{k,j} = \bgamma_{k+1,j}$.
\ElsIf {$v_j, v_{j_1}\in V_{k}$ collided at $t_k$}
\State Perform backward collision between $\bgamma_{k+1,j}$ and $\bgamma_{k+1,j_1}$ and obtain $\bgamma_{k,j}$ and $\bgamma_{k,j_1}$ following~\eqref{gammaEqtn} and~\eqref{gammaTildeEqtn}. 
\EndIf
\State Obtain $\{\bgamma_{k,1},\bgamma_{k,2}, \ldots, \bgamma_{k,N}\}$.
\EndFor
\end{algorithmic}
\end{algorithm}

Like the continuous adjoint system~(\Cref{ContinuousFormulation}), by solving the direct equations and the adjoint equations only once, we obtain all components of the gradient $\partial_\alpha J$, for any value of the dimensionality of $\alpha$. We remark that the choice of collision times, collision partners, and collision angles is not included in the variational principle. These parts of the DSMC method are determined externally.


The restriction to Maxwell molecules significantly simplifies the adjoint DSMC method presented in this section. Since the collision kernel $q$ is a constant, the collision rate for a pair of velocities $v$ and $v_1$ and the choice of the collision parameter $\sigma$ do not depend on the values of the velocities. This implies that variations  $\delta v$ and  $\delta v_1$ in the velocities do not cause variations in the choice of collision pairs nor in the collision parameter, which simplifies the adjoint DSMC equation~\eqref{gammaEqtns}. This is the main reason that we restricted our attention to Maxwell molecules in this work. 
We believe it is possible to extend the adjoint DSMC method to non-Maxwell interactions. Extension to collision kernels $q=q(\theta)$ with angular dependence should be straightforward since, in that case, the only thing that changes is the sampling of the collision parameter $\sigma$, while it is still independent of pre-collision velocities $v$ and $v_1$. Extensions to the VHS model with $q=q(|v-v_1|)$ collision kernel and other non-Maxwell interactions are more intricate since either the collision rates or the collision parameter $\sigma$ will depend on pre-collision velocities. The former could be mitigated by using a different DSMC method with a constant collision rate such as \cite{BobylevNanbu2000}, while the latter would produce additional terms in $\partial_{w_k}w_k', \partial_{\tilde w_k}w'_k, \partial_{w_k}\tilde w_k', \partial_{\tilde w_k}\tilde w'_k$ in~\eqref{eq:fwd_perturb45}. We leave such generalizations for future work and they are beyond the scope of this paper.






\section{Relationship between the continuous adjoint and the adjoint DSMC formulation}
\label{sec:discussion}
In this section, we discuss the similarities and differences between the continuous adjoint formulation (\Cref{ContinuousFormulation}) and the adjoint DSMC formulation (\Cref{DiscreteFormulation}). The former is derived under the ``optimize-then-discretize'' (OTD) approach, while the latter is based on the ``discretize-then-optimize'' (DTO) approach. For a given discretization scheme, OTD and DTO may not be equivalent, as demonstrated by many examples in the literature~\cite{hager2000runge,burkardt2002insensitive,
abraham2004effect,hinze2012discretization,liu2019non}. The DTO approach is particularly preferred for problems that are otherwise unsolvable under the OTD approach, such as the heat-transfer optimization problem~\cite{betts2005discretize,ghobadi2009discretize}.

Since DTO and OTD are generally not the same, we investigate the analytical properties of the two adjoint systems. The analysis in~\Cref{sec:true_continuous} illustrates the essential role of the continuous adjoint variable $\gamma (v,t)$, as a Fr\'echet derivative of the objective function. We find a similar result for  the discrete adjoint variable $\bgamma_{k,i}$ for DSMC  in~\Cref{sec:DSMC_property}. Based on these results, we find a direct relationship between $\gamma (v,t)$ and $\bgamma_{k,i}$ in~\Cref{sec:connections} that connects the two adjoint systems.

We remark that the derivatives in this section have slightly different meanings from the ones used before. In~\Cref{ContinuousFormulation} and~\Cref{DiscreteFormulation}, the constraints for the state variables are not directly applied, but instead are imposed by the adjoint variables in the Lagrangian formulation. Here in~\Cref{sec:discussion}, the derivatives with respect to $f(v)$ or $v_k$ are computed with the assumption that they are directly constrained to satisfy the Boltzmann equation~\eqref{eq:homoBoltz} or the DSMC equations, respectively.

\subsection{Continuous adjoint variable as a Fr\'echet derivative}\label{sec:true_continuous}

If $f(v,t)$ is a solution of the  Boltzmann equation~\eqref{eq:homoBoltz}, then the linearized equation for a perturbed distribution $f(v,t) + \delta f(v,t)$ is
\begin{equation}\label{eq:new_linear2}
    \partial_t \delta f(v,t) = f(v,t)\, L[f(v,t)]\left(\frac{\delta f(v,t)}{f(v,t)}\right)
\end{equation}
using  the linearized  collision operator~\eqref{eq:linearQ}.

For objective function $J_1(\alpha) = \int_v r(v) f(v,T)dv$, the perturbation $\delta f$ causes a perturbation  $\delta J_1$, which satisfies 
$$\delta J_1 =  \int_\Rthree \delta f(v,t) \frac{\delta J_1}{\delta f(v,t)} dv$$
for any $t$ because the perturbation $\delta f(t, \cdot)$  determines $\delta f (t', \cdot)$ for any $t'>t$.
Since $J_1$ is time-independent, $\forall t\in[0,T]$, we have
\begin{align*}
    0 =\partial_t \delta J_1 &= \partial_t \left( \int_\Rthree \delta f(v,t) \frac{\delta J_1}{\delta f(v,t)} dv \right) \\
    &=  \int_\Rthree (\partial_t \delta f)  \left(\frac{\delta J_1}{\delta f}\right)dv + \int_\Rthree \delta f\, \partial_t \left(\frac{\delta J_1}{\delta f}\right)dv\\
    &=\int_\Rthree f\, L[f]\left(\frac{\delta f}{f}\right) \,\frac{\delta J_1}{\delta f}dv + \int_\Rthree \delta f \partial_t \left(\frac{\delta J_1}{\delta f}\right)dv\\
     &=\int_\Rthree f\, L^*[f] \left(\frac{\delta J_1}{\delta f} \right) \,\frac{\delta f}{f} dv + \int_\Rthree \delta f \partial_t \left(\frac{\delta J_1}{\delta f}\right)dv\\
     &=\int_\Rthree \delta f\left( L^*[f] \left(\frac{\delta J_1}{\delta f} \right) +  \partial_t \left(\frac{\delta J_1}{\delta f}\right)\right)dv.
\end{align*}
Here, we have used the chain rule, the duality between $L[f]$ and $L^*[f]$, and~\eqref{eq:new_linear2}.

Since the perturbation $\delta f$ is arbitrary, the following holds for any $t$
\begin{equation}\label{eq:new_adjoint}
    - \partial_t \left(- \frac{\delta J_1}{\delta f}\right) = L^*[f] \left(-\frac{\delta J_1}{\delta f} \right).
\end{equation}
Note that equation~\eqref{eq:new_adjoint}  for $-\frac{\delta J_1}{\delta f}$  is the same as equation~\eqref{backwards1} for  the continuous adjoint  function $\gamma(v,t)$. Additionally, as seen directly from the objective function, the ``final'' condition for~\eqref{eq:new_adjoint} is $- \frac{\delta J_1}{\delta f(v,T)} = - r(v)$, the same as the ``final'' condition~\eqref{backwards0} for $\gamma(v,T)$. 
Based on these two facts, we obtain the main result of this subsection:
\begin{equation}\label{eq:true_continuous_gamma}
    \gamma(v,t) = - \frac{\delta J_1}{\delta f(v,t)},\quad t\in[0,T].
\end{equation}
Therefore, the continuous adjoint variable $\gamma(v,t)$ is the negative Fr\'echet derivative of the objective function with respect to the state variable $f(v,t)$.



\subsection{DSMC adjoint variable as a derivative}~\label{sec:DSMC_property}
The analysis of the DSMC adjoint variable $\bgamma$ is a discrete version of the analysis in the previous section for the continuous adjoint variable $\gamma$. Consider  small perturbations $\delta w_k$ and $\delta \tilde w_k$ in the pre-collision velocities $w_k$ and $\tilde w_k$. The resulting first-order variations in the post-collision velocities are $\delta w_{k+1}$ and $\delta \tilde w_{k+1}$. Since no other velocities are changed at times $t_k$ and $t_{k+1}$ and because $\cJ_1$ is time independent, then
\begin{equation*} 
    \delta \cJ_1 = \partial_{w_k} \cJ_1 \cdot \delta w_k + \partial_{\tilde w_k} \cJ_1  \cdot \delta \tilde w_k = \partial_{w_{k+1}} \cJ_1 \cdot \delta w_{k+1} + \partial_{\tilde w_{k+1}} \cJ_1 \cdot \delta\tilde w_{k+1},
\end{equation*}
which can be rewritten in operator form, and then by using~\eqref{eq:fwd_perturb_General}, as
\begin{equation*} %
        \begin{pmatrix}
   \partial_{w_k} \cJ_1 ^T &
    \partial_{\tilde w_k} \cJ_1 ^T
    \end{pmatrix}  
         \begin{pmatrix}
   \delta{w_k} \\
    \delta{\tilde w_k}  
    \end{pmatrix}  
    \nonumber 
   =
        \begin{pmatrix}
   \partial_{w_{k+1}} \cJ_1 ^T &
    \partial_{\tilde w_{k+1}} \cJ_1 ^T
    \end{pmatrix}  
         \begin{pmatrix}
   \delta{w_{k+1}} \\
    \delta{\tilde w_{k+1}}  
    \end{pmatrix}  
     =
        \begin{pmatrix}
   \partial_{w_{k+1}} \cJ_1 ^T &
    \partial_{\tilde w_{k+1}} \cJ_1 ^T
    \end{pmatrix}  
   A(\sigma_k,\alpha_k) 
         \begin{pmatrix}
   \delta{w_k} \\
    \delta{\tilde w_k}  
    \end{pmatrix}  .
    \end{equation*}
  Since this is true for any values of $\delta w_k$ and $\delta \tilde w_k$, it follows that 
\begin{equation*} %
        \begin{pmatrix}
   \partial_{w_k} \cJ_1 ^T &
    \partial_{\tilde w_k} \cJ_1 ^T
    \end{pmatrix}  
=          \begin{pmatrix}
   \partial_{w_{k+1}} \cJ_1 ^T &
    \partial_{\tilde w_{k+1}} \cJ_1 ^T
    \end{pmatrix}  
   A(\sigma_k,\alpha_k) .
    \end{equation*}
By taking the transpose and using $A^T=B$, then
\begin{equation}\label{eq:fwd_perturb3}
        \begin{pmatrix}
   \partial_{w_k} \cJ_1\\
    \partial_{\tilde w_k} \cJ_1 
    \end{pmatrix}  = B(\sigma_k,\alpha_k) 
    \begin{pmatrix}
    \partial_{w_{k+1}} \cJ_1 \\
    \partial_{\tilde w_{k+1}} \cJ_1
    \end{pmatrix}.
\end{equation}
Equation~\eqref{eq:fwd_perturb3} shows that $(\partial_{w_k} \cJ_1,\partial_{\tilde w_k} \cJ_1)$ satisfies the same back-propagation rule as~\eqref{eq:AB_gamma} for 
$(\bgamma_k, {\tilde \bgamma_k})$.

Recall~\eqref{gammaFinal} for the final condition for the adjoint DSMC particles, 
\begin{equation}\label{eq:fwd_perturb7}
\bgamma_F = -\partial_v r(w_F)  =  - N\partial_{w_F} \cJ_1.
\end{equation}
in which the second equation comes from the definition of $\cJ_1$ in~\eqref{eq:DTO_obj}.
With the same final condition~\eqref{eq:fwd_perturb7} (up to a constant scaling $-N$), and the same back-propagating rule~\eqref{eq:AB_gamma} and~\eqref{eq:fwd_perturb3} from $t_{k+1}$ to $t_k$, we conclude that \begin{equation}~\label{eq:true_DSMC_gamma}
    \bgamma_{k,i} = -N\frac{\partial \cJ_1}{\partial v_{k,i}},\quad \forall k,i,\quad\text{where } \cJ_1 = \frac{1}{N}\sum_{i=1}^N r(v_{F,i}).
\end{equation}

\subsection{Connections between the continuous and DSMC adjoint variables}\label{sec:connections}
Given an empirical distribution and its limit, we know by the strong law of large numbers that
\begin{equation*}
\frac{1}{N}\sum_{i=1}^N I_\Omega(v_{F,i}) \xrightarrow[]{a.s.} \int_\Omega f(v,T) dv
\end{equation*}
for every measurable set $\Omega\subseteq \Rthree$, where $I_\Omega$ is the indicator function. Also,
\begin{equation*}
\lim_{N\rightarrow \infty}\cJ_1 =\lim_{N\rightarrow \infty}\frac{1}{N}\sum_{i=1}^N r(v_{F,i}) = \int_{\Rthree} r(v)f(v,T)dv = J_1.
\end{equation*}
Next, we derive the connections between the two derivatives at any time step $t_k$,
$$
\frac{\delta J_1}{\delta f(v,t_k)} \quad\text{and}\quad \frac{\partial \cJ_1}{\partial v_{k,i}},
$$
which will directly uncover the relationship between the continuous and DSMC adjoint variables based on~\eqref{eq:true_continuous_gamma} and~\eqref{eq:true_DSMC_gamma}.

We consider two empirical density functions,
$$ \frac{1}{N} \sum_{i=1}^N \delta(v-v_{k,i})\quad\text{and}\quad \frac{1}{N} \sum_{i=1}^N \delta(v-\bar v_{k,i}),\quad \bar v_{k,i} = v_{k,i} + \delta v_{k,i},$$ 
as the approximations to the density functions $f(v,t_k)$ and $f(v,t_k) +\delta f(v,t_k)$, respectively. We remark that $v_{k,i}$ are random variables following the distribution function $f(v,t_k)$, while we choose $\delta v_{k,i}$ to be  deterministic, as follows:
\begin{equation*}
    \delta v_{k,i} = \begin{cases}
    \eta, & v_{k,i}\in \Omega,\\
    0, & v_{k,i}\in \Rthree \backslash \Omega.
    \end{cases}
\end{equation*}
Here, $\Omega$ is an arbitrary measurable set and $\eta$ is a small constant vector. 

The first-order variation of the continuous and the discrete objective functions can be stated as below. For any time step $t_k$,
\begin{equation} \label{eq:dJ0}
    \delta J_1 = \int_\Rthree \frac{\delta J_1}{\delta f}  \delta f dv= \int_\Rthree \frac{\delta J_1}{\delta f(v,t_k)} \bigg( f(v,t_k) +\delta f(v,t_k) - f(v,t_k) \bigg) dv,
\end{equation}
\begin{equation}\label{eq:dJ1}
\delta \cJ_1 =  \sum_{i=1}^N \delta v_{k,i} \cdot \frac{\partial \cJ_1}{\partial v_{k,i}}.
\end{equation}

Consider $v_i = v_{k,i}$ as one of the $N$ random variables following the distribution $f(v,t_k)$, and denote $\phi(v) =   \frac{\delta J_1}{\delta f(v,t_k)}$. Based on~\eqref{eq:dJ0}, we have
\begin{align}
\delta J_1  &=\int_\Rthree \frac{\delta J_1}{\delta f(v,t_k)} \big( f(v,t_k) +\delta f(v,t_k) \big) dv - \int_\Rthree \frac{\delta J_1}{\delta f(v,t_k)}  f(v,t_k) dv, \nonumber\\
&=  \mathbb E_{\bar{v}_i} \bigg[\frac{\delta J_1}{\delta f(\bar{v}_i,t_k)}\bigg] -  \mathbb E_{v_i} \bigg[\frac{\delta J_1}{\delta f(v_i,t_k)}\bigg] 
\nonumber  \\
&=  \mathbb E_{\bar{v}_i} [\phi(\bar{v}_i)] -  \mathbb E_{v_i} [\phi({v}_i)]
\nonumber\\
&=  \mathbb E_{ v_i } \big [ I_\Omega(v_i)  \phi( v_i + \eta) + I_{\Rthree\backslash \Omega}(v_i) \phi( v_i) \big] -\mathbb E_{v_i} [\phi({v}_i)]
\nonumber \\
&=  \mathbb E_{ v_i } [ I_\Omega(v_i)  \phi( v_i + \eta) ] -\mathbb E_{v_i} [I_\Omega(v_i) \phi({v}_i)],\nonumber\\
&=  \mathbb E_{ v_i } [ I_\Omega(v_i)  \left(\phi( v_i + \eta) -\phi(v_i)\right)] \nonumber \\
&\approx   \eta \cdot \mathbb E_{ v_i } \big[ I_\Omega(v_i)  \phi'(v_i) \big]
\label{eq:dJ2}
\end{align}
On the other hand,~\eqref{eq:dJ0} is related to~\eqref{eq:dJ1} as the following
\begin{equation}\label{eq:dJ3}
    \delta J_1 = \mathbb E_{v_1,\dots,v_N}[\delta \cJ_1] =  \sum_{i=1}^N  \mathbb E_{v_i }\bigg[\delta v_i \cdot \frac{\partial \cJ_1}{\partial v_i}\bigg] = N \eta \cdot   \mathbb E_{v_i }\bigg[I_\Omega(v_i) \frac{\partial \cJ_1}{\partial v_i}\bigg].
\end{equation}
We assume that all $v_i = v_{k,i}$, for $i = 1,\dots,N$, are (approximately) i.i.d. random variables following the same velocity distribution $f(v,t_k)$.

Therefore, by combining~\eqref{eq:dJ2} and~\eqref{eq:dJ3}, we have
\[
\eta \cdot \mathbb E_{ v_i } \big[ I_\Omega(v_i)  \phi'(v_i) \big] = N \eta \cdot   \mathbb E_{v_i }\bigg[I_\Omega(v_i) \frac{\partial \cJ_1}{\partial v_i}\bigg].
\]
Since $\eta$ is an arbitrarily small vector in $\Rthree$, then 
\begin{equation*}
  \mathbb E_v \bigg[ I_\Omega(v)\phi'(v) \bigg] =N \mathbb E_{v_i }\bigg[ I_\Omega(v_i) \frac{\partial \cJ_1}{\partial v_i}\bigg]
  = N\mathbb E_v  \bigg[  \mathbb E_{v_i} \bigg[I_\Omega(v_i) \frac{\partial \cJ_1}{\partial v_i} \bigg| v_i = v \bigg] \bigg] = \mathbb E_v \bigg[ I_\Omega(v) \mathbb E_{v_i} \bigg[N \frac{\partial \cJ_1}{\partial v_i} \bigg| v_i = v \bigg] \bigg],
\end{equation*}
where we change the variable from $v_i$ to $v$ for the first term and apply the law of total expectation to the second term.

Since the set $\Omega$ is arbitrary, $\phi'(v) =(\frac{\delta J_1}{\delta f(v,t_k)})' = -\gamma'(v,t_k)$ and $ N\frac{\partial \cJ_1}{\partial v_{k,i}} =- \bgamma_{k,i}$, which have been shown previously in~\eqref{eq:true_continuous_gamma} and~\eqref{eq:true_DSMC_gamma}, we obtain the following equation 
\begin{equation} \label{eq:dJ4}
\gamma'(v,t_k)   = \mathbb E[\bgamma_{k,i}|v_{k,i}=v].
\end{equation}
\Cref{eq:dJ4} serves as a bridge connecting the two adjoint systems that we have derived in~\Cref{ContinuousFormulation} and~\Cref{DiscreteFormulation}.

\subsection{The unbiasedness of the adjoint DSMC gradient estimator}
In~\eqref{optimization1} and~\eqref{DalphaEqtn}, we derive two gradient formulae, $\partial_\alpha J$ and $\partial_\alpha \cJ$, through the continuous and the discrete adjoint systems, respectively. Next, we show the unbiasedness of the gradient estimator $\partial_\alpha \cJ$ computed through the adjoint DSMC approach:
\begin{equation}
    \partial_\alpha J = \mathbb{E}_{v_{0,1},\ldots, v_{0,N}\sim f_0(v;\alpha)}[\partial_\alpha \cJ]. 
\end{equation}

There are a direct way and an indirect way to draw samples from continuous distributions $f(v,0) = f_0(v; \alpha)$, which is also our initial distribution for the Boltzmann equation:
\[
\text{Direct way: }v_{0,i} \sim f_0(v;\alpha)\]
\[
\text{Indirect way: }v_{0,i} = g(\epsilon_{0,i},\alpha),\, \epsilon_{0,i} \sim \phi(\epsilon).
\]
The second and indirect approach is to first sample from a simpler base distribution $\phi(\epsilon)$, which is independent of the parameter $\alpha$, and then transform this variate through a deterministic path $ g(\epsilon,\alpha)$. This is often referred to as the sampling path or the sampling process~\cite{glasserman2013monte,mohamed2019monte}. For an invertible path, we have the mass-preserving equation:
\begin{equation} \label{eq:mass-preserving}
   \phi(\epsilon) =   f_0\Big(g(\epsilon,\alpha);\alpha\Big) \det{ \Big(\nabla_\epsilon g(\epsilon,\alpha) \Big)}  .
\end{equation} 

Recall in~\Cref{ContinuousFormulation} where we derive the OTD framework, 
\[\partial_\alpha J = \partial_\alpha J_2 = -\partial_\alpha \int_\Rthree \gamma(v,0)f_0(v;\alpha) dv = -\partial_\alpha \mathbb{E}_{v\sim f_0(v;\alpha)}[\gamma(v,0)].
\]
Based on the Law of the Unconscious Statistician (LOTUS),
\begin{equation}
    \mathbb{E}_{v\sim f_0(v;\alpha)}[\gamma(v,0)] = \mathbb{E}_{\epsilon\sim \phi(\epsilon)}[\gamma(g(\epsilon,\alpha),0)].
\end{equation}
Therefore, we have
\begin{align*}
    \partial_\alpha J  &=  -\partial_\alpha \mathbb{E}_{v\sim f_0(v;\alpha)}[\gamma(v,0)]\\
    & =  -\partial_\alpha \mathbb{E}_{\epsilon\sim \phi(\epsilon)}\bigg[\gamma\Big(g(\epsilon,\alpha),0\Big)\bigg]= -\partial_\alpha \left(\int \phi(\epsilon) \gamma \Big(g(\epsilon,\alpha),0 \Big) d\epsilon\right)\\
    & = -\int \phi(\epsilon)\left( \frac{\partial \gamma(v,0)}{\partial v}\bigg|_{v=g(\epsilon,\alpha)}\cdot \frac{\partial g(\epsilon,\alpha)}{\partial \alpha}\right) d\epsilon\\
     & = -\int f_0\Big(g(\epsilon,\alpha); \alpha\Big) \left( \frac{\partial \gamma(v,0)}{\partial v}\bigg|_{v=g(\epsilon,\alpha)}\cdot \frac{\partial g(\epsilon,\alpha)}{\partial \alpha}\right) \det{ \Big(\nabla_\epsilon g(\epsilon,\alpha) \Big)} d\epsilon \quad \text{\big(by ~\Cref{eq:mass-preserving}\big)}\\
    & = -\int f_0(v;\alpha)\left( \frac{\partial \gamma(v,0)}{\partial v} \cdot \frac{\partial v}{\partial \alpha}\right) dv.
\end{align*}
The last equality is based on the change of variable $v = g(\epsilon,\alpha)$. Together with~\Cref{eq:dJ4}, we have
\begin{align*}
    \partial_\alpha J  & = -\int f_0(v;\alpha) \left(\frac{\partial \gamma(v,0)}{\partial v} \cdot \frac{\partial v}{\partial \alpha} \right)dv  = - \frac{1}{N}\sum_{i=1}^N \mathbb {E}_{v_{0,i}\sim f_0(v;\alpha)} \bigg[ \bgamma_{0,i} \cdot \frac{\partial v_{0,i}}{\partial \alpha} \bigg]\\
& =  \mathbb{E}_{v_{0,1},\ldots, v_{0,N}\sim f_0(v;\alpha)} \bigg[- \frac{1}{N}\sum_{i=1}^N \bgamma_{0,i}\cdot \frac{\partial v_{0,i}}{\partial \alpha} \bigg]=  \mathbb{E}_{v_{0,1},\ldots, v_{0,N}\sim f_0(v;\alpha)} \, [\partial_\alpha \cJ].
\end{align*}
In conclusion, $\partial_\alpha \cJ$ is an unbiased pathwise Monte Carlo gradient estimator for the objective function $J$.

If we solve the continuous adjoint equation~\eqref{backwards1} numerically and then compute the gradient for $\alpha$ following~\eqref{optimization1}, the result should match the gradient computed by the DSMC adjoint approach in~\Cref{DiscreteFormulation} within numerical error. For example,~\eqref{backwards1} can be solved numerically by the DSMC-type scheme which we propose in~\Cref{sec:DSMCtype_method}. The continuous and the DSMC adjoint systems are equivalent in optimizing the model parameter $\alpha$. We will present several numerical examples that verify the gradient accuracy and demonstrate the optimization process in~\Cref{sec:numerics}.



\section{Numerical results}\label{sec:numerics}
In this section we discuss the results of numerical simulations computing the gradients of the objective function \eqref{eq:DTO_obj}, $\cJ_1 =\frac{1}{N}\sum_{i=1}^N r(v_{F,i}) \approx \int_{\bbR^3} r(v) f(v,T) dv = J_1$, at the final time $t=T$ with respect to the parameter $\alpha$ in the initial conditions $f_0(v;\alpha)$. Four different methods of the gradient computation are used here: (i) finite difference method using several forward DSMC simulations with different parameter values, (ii) the adjoint DSMC method, (iii) the DSMC-type scheme for the adjoint equation, (iv) the direct discretization of the continuous adjoint equation~\eqref{backwards1}. All four methods lead to the same gradient values, but the adjoint DSMC method is the best in terms of performance given we do not need many digits of accuracy.

Here we consider Maxwellian gas, so the distribution function $f(v,t)$ obeys the Boltzmann equation \eqref{eq:homoBoltz} with a collision operator kernel $q(v-v_1,\sigma) = q(\sigma)$. We further assume that $q(\sigma)=1/(4\pi)$ and $\rho(t)=\int_{\bbR^3} f(v,t) dv=1$ (conserved throughout the evolution of the distribution function), and thus $\mu =\rho \int_\Stwo q(\sigma)d\sigma = 1$ in~\eqref{eq:homoBoltz_DSMC_1}.

For the function $r(v)$ in Section~\ref{sec:Forward_DSMC}-\ref{sec:Direct_num}, we use $v_l^2$ and $v_l^4$, $l\in\{x,y,z\},$ so the objective functions are 
$$m2_l(t_k) \triangleq T_l(t_k)=\frac{1}{N}\sum_{i=1}^N{(v_{k,i}^l)^2 } \approx \int_{\bbR^3} v_l^2 f(v,t_k) dv,\quad l\in\{x,y,z\},$$ $$m4_l(t_k) \triangleq\frac{1}{N}\sum_{i=1}^N{(v_{k,i}^l)^4 } \approx \int_{\bbR^3} v_l^4 f(v,t_k) dv, \quad l\in\{x,y,z\}.$$
They are the second-order and the fourth-order velocity moments of the distribution function in the $l$-direction at the time $t=t_k$. We have six objective functions in total.  For the parameter $\alpha$ we use temperature values in the initial distribution function $\alpha=[T_x^0, T_y^0, T_z^0]$. We further refer to these gradients as $\frac{\partial T_l}{\partial T_p^0}$ and $\frac{\partial m4_l}{\partial T_p^0}$ respectively, $l,p\in\{x,y,z\}$. 
Here the dimension of $\alpha$ is only 3. Still, a real advantage of these methods is, of course, when the vector $\alpha$ is highly multi-dimensional since the described methods allow one to compute all the components of the gradient $\frac{\delta \cJ_1}{\delta \alpha}$ by doing the only \textit{one} forward DSMC simulation and \textit{one} backward adjoint DSMC simulation.

For all the methods we use the same initial condition, anisotropic Gaussian, 
\begin{equation}~\label{eq:IC}
f_0(v) = \frac{1}{(2\pi)^{3/2}\sqrt{T_x^0 T_y^0 T_z^0}} \exp\left(-\frac{v_x^2}{2T_x^0}-\frac{v_y^2}{2T_y^0}-\frac{v_z^2}{2T_z^0}\right),
\end{equation}
where $T_x^0=0.5,T_y^0=1,T_z^0=1$, and the total density $\rho=1$. In this case, the solution to the Boltzmann equation~\eqref{eq:homoBoltz} will relax to an isotropic Gaussian with the temperature $T_M=(T_x^0+T_y^0+T_z^0)/3=0.8333(3)$ over time. In all the tests below, we use the forward Euler time-integration scheme with a time-step $\Delta t=0.1$. The gradients, $\frac{\partial T_l}{\partial T_p^0}$ and  $\frac{\partial m4_l}{\partial T_p^0}$, are computed at the final time $t=T=2$.


\subsection{Forward DSMC simulation} \label{sec:Forward_DSMC}
\begin{figure}
\centering
\includegraphics[width=0.4\textwidth]{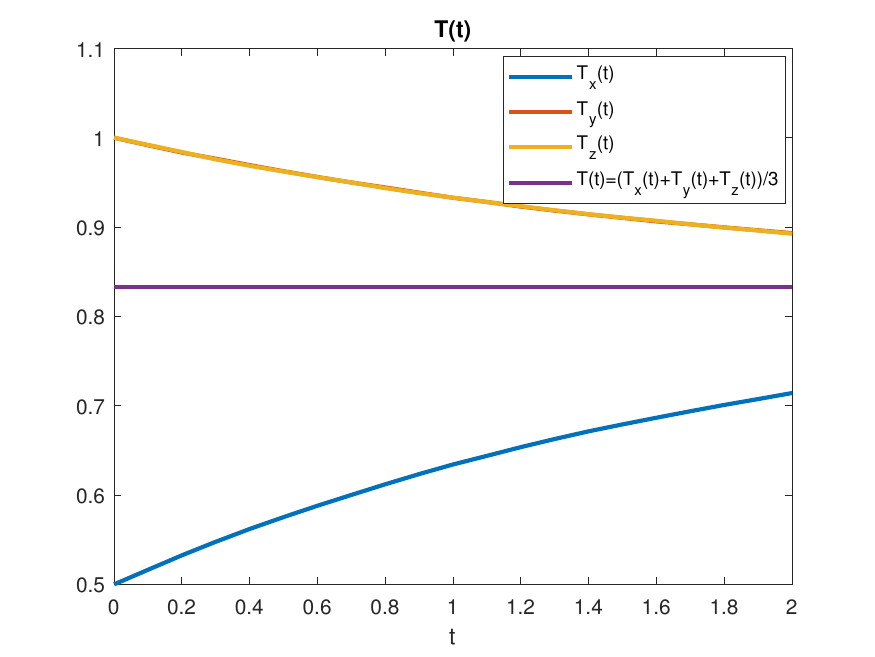}
\includegraphics[width=0.4\textwidth]{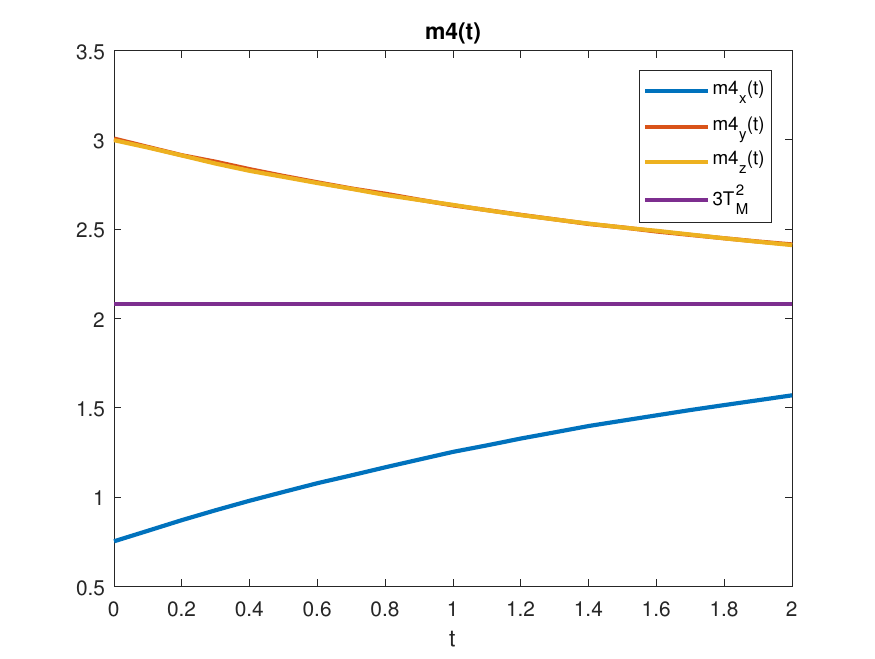}
\caption{The relaxation of temperatures $T_x(t),T_y(t),T_z(t)$ to the value $T_M=(T_x^0+T_y^0+T_z^0)/3=0.8333(3)$ as well as the average temperature $(T_x(t)+T_y(t)+T_z(t))/3$ that stays constant throughout the simulation (left) and the relaxation of fourth moments $m4_x(t),m4_y(t),m4_z(t)$ to the value $3T_M^2=2.08333(3)$ from DSMC simulation with initial condition \eqref{eq:IC}, $N=10^6$ particles and $\Delta t=0.1$. The red and yellow lines on both graphs practically coincide.} \label{fig:Forward_DSMC}
\end{figure}

We solve the Boltzmann equation~\eqref{eq:homoBoltz} with the Nanbu--Babovsky method as described in~\Cref{alg:DSMC}. We represent the distribution function with $N$ particles,  ranging from $10^6$ to $10^8$, and sample the initial condition from the distribution in~\eqref{eq:IC}. Collisions are performed according to the  collision rules~\eqref{BoltzmannSolution} and the collision angles are sampled uniformly over a unit sphere according to the collision kernel of Maxwellian particles, $q(v-v_1,\sigma) =1/(4\pi)$. Based on~\Cref{alg:DSMC}, the fraction of particles that collide at every time step is $N_c/N = \Delta t \mu$ and is equal to $10\%$ for $\Delta t=0.1$, $\mu=1$.
The total kinetic energy $K(t)=T_x(t)+T_y(t)+T_z(t)$ 
and the total momentum $p(t)=(p_x(t),p_y(t),p_z(t))$, where  $p_l(t)=\int_{\bbR^3} v_l f(v,t) dv$,
are conserved by construction of the algorithm, since every pair-wise collision is elastic.
Fig.~\ref{fig:Forward_DSMC} shows the relaxation of temperatures $T_x(t),T_y(t),T_z(t)$ towards the value $T_M$ as well as the relaxation of the fourth-order moments $m4_x(t),m4_y(t),m4_z(t)$ towards the value $3T_M^2$ as functions of time.

Running the standard DSMC simulations several times with randomly sampled initial conditions allows us to estimate the mean values and the standard deviations of quantities of interest at $t=T=2$. Details are shown in \Cref{Table1}, where the standard deviation $\sigma_X$ is computed using $M_s=100$ simulations with $N=10^6, 10^8$ particles in each.
The standard deviations of the quantities of interest give us estimates of the random errors in DSMC algorithm with $N$ particles. They scale with $1/\sqrt{N}$ which can be seen from the $\sigma_X$ values in~\Cref{Table1}. The actual standard deviations of the mean values $\bar{X}$ in~\Cref{Table1} computed using $M_s$ independent simulations are approximately $1/\sqrt{M_s}$ times the values of $\sigma_X$ for one simulation from~\Cref{Table1}. That is, $\sigma_{\bar{X}} = {\sigma_X}/  {\sqrt{M_s}}$. We further estimate random errors in the mean values $\bar{X}$ by using $95\%$ of the trust interval of radius $2\sigma_{\bar{X}}$, so $e^{rand}_{\bar{X}}=2\sigma_{\bar{X}}=\frac{2 \sigma_X}{\sqrt{M_s}}$. Finally, we have an estimates of expectations:
\begin{equation}~\label{eq:E[X]approx}
\mathbb E [X] \approx  \bar{X}  \pm \frac{2 \sigma_X}{\sqrt{M_s}}. 
\end{equation}
Based on this approach and the values of $\sigma_X$ from \Cref{Table1}, we can see that when $N=10^8$ and $M_s=100$, the random errors of $T_x(t=2), T_y(t=2), T_z(t=2)$ are $e^{rand}_{T_l} = 2 \sigma_{T_l}/\sqrt{M_s} \approx 0.00002$ (or relative error $e^{rand}_{T_l}/T_l \approx 0.002\% $). The random errors of $m4_x(t=2), m4_y(t=2), m4_z(t=2)$ are $e^{rand}_{m4_l} =2 \sigma_{m4_l}/\sqrt{M_s} \approx 0.0001$ (or relative error $e^{rand}_{m4_l}/m4_l \approx 0.005\% $), $l\in\{x,y,z\}$.

 \begin{table} 
\small \renewcommand{\arraystretch}{1.0}
\centering 
\begin{tabular}{|c|c|c|c|c|}  
  \hline
  &  \multicolumn{4}{|c|}{$M_s=100$} \\ 
  \cline{2-5}
  &  \multicolumn{2}{|c|}{$N=10^6$} &  \multicolumn{2}{|c|}{$N=10^8$} \\ 
  \cline{2-5}
  & \rule{0pt}{10pt}  $\bar{X}$ & $\sigma_{X}$ & $\bar{X}$ & $\sigma_{X}$ \\
  \hline
$p_x(t=2)$  & 3.5237e-18 & 1.7446e-17  & 2.6791e-19 & 1.9801e-17 \\
$p_y(t=2)$  & 8.3766e-19 & 2.4584e-17  & 1.4893e-18 & 2.9098e-17  \\
$p_z(t=2)$  &-3.5049e-18 & 2.6808e-17  & 4.8894e-18 & 2.8732e-17  \\ 
$T_x(t=2)$  & 0.71405 & 0.00077  & 0.7138388 & 0.0000735 \\
$T_y(t=2)$  & 0.89289 & 0.00087  & 0.8930843 & 0.0000890 \\
$T_z(t=2)$  & 0.89306 & 0.00105  & 0.8930770 & 0.0001009 \\
$m4_x(t=2)$ & 1.5699  & 0.0043   & 1.5690043 & 0.0004608 \\ 
$m4_y(t=2)$ & 2.4097  & 0.0059   & 2.4114132 & 0.0005268 \\  
$m4_z(t=2)$ & 2.4115  & 0.0070   & 2.4113140 & 0.0006946 \\
   \hline
\end{tabular}
  \caption{Mean values of $p_x, p_y, p_z$, $T_x, T_y, T_z$ and $m4_x, m4_y, m4_z$ and their standard deviations at $t=2$, which are computed using $M_s$ standard DSMC simulations with $N$ particles in each. The temperatures in the initial condition~\eqref{eq:IC} are set to be $[T_x^0 \; T_y^0 \; T_z^0]=[0.5 \; 1 \; 1]$.}\label{Table1}
\end{table}
 
To directly compute the gradients of the velocity moments with respect to the parameter $\alpha=[T_x^0, T_y^0, T_z^0]$ by finite difference, we need to perturb the parameter by a small amount and then compute the corresponding target values at those different initial temperatures. For each component in $[T_x^0 \; T_y^0 \; T_z^0]$, we apply a perturbation of size $\Delta T_l^0=0.1$ or $0.05$ away from the original values of $[T_x^0 \; T_y^0 \; T_z^0]=[0.5 \; 1 \; 1]$. In~\Cref{Table2}, we gather the resulting mean values of the moments that are computed using $M_s=100$ simulations and $N=10^8$ particles in each simulation.

\begin{table} 
\small 
\setlength{\tabcolsep}{2pt}
\centering 
\begin{tabular}{|c|c|c|c|c|c|c|c|} 
  \hline
$[T_x^0 \; T_y^0 \; T_z^0]$ & [0.5 1 1] & [0.4 1 1] & [0.6 1 1] & [0.5 0.9 1] & [0.5 1.1 1] & [0.5 1 0.9] & [0.5 1 1.1] \\
  \hline
$T_x(t=2)$  &  0.71384 &  0.65661 &  0.77108 &  0.69246 &  0.73523 &  0.69246 &  0.73523 \\
$T_y(t=2)$  &  0.89308 &  0.87171 &  0.91447 &  0.83586 &  0.95032 &  0.87171 &  0.91447 \\
$T_z(t=2)$  &  0.89308 &  0.87169 &  0.91445 &  0.87169 &  0.91445 &  0.83584 &  0.95030 \\
$m4_x(t=2)$ &  1.56900 &  1.35147 &  1.80946 &  1.47149 &  1.67080 &  1.47149 &  1.67080 \\
$m4_y(t=2)$ &  2.41141 &  2.30639 &  2.52069 &  2.10811 &  2.73762 &  2.29959 &  2.52749 \\
$m4_z(t=2)$ &  2.41131 &  2.30628 &  2.52058 &  2.29949 &  2.52737 &  2.10801 &  2.73750 \\
  \hline
$[T_x^0 \; T_y^0 \; T_z^0]$ &   & [0.45 1 1] & [0.55 1 1] & [0.5 0.95 1] & [0.5 1.05 1] & [0.5 1 0.95] & [0.5 1 1.05] \\
  \hline
$T_x(t=2)$  &   & 0.68523 & 0.74247 & 0.70315 & 0.72453 & 0.70315 & 0.72453 \\
$T_y(t=2)$  &   & 0.88240 & 0.90376 & 0.86447 & 0.92171 & 0.88240 & 0.90378 \\
$T_z(t=2)$  &   & 0.88238 & 0.90378 & 0.88238 & 0.90376 & 0.86445 & 0.92168 \\
$m4_x(t=2)$ &   & 1.45739 & 1.68646 & 1.51973 & 1.61938 & 1.51973 & 1.61938 \\
$m4_y(t=2)$ &   & 2.35837 & 2.46534 & 2.25690 & 2.57166 & 2.35497 & 2.46893 \\
$m4_z(t=2)$ &   & 2.35827 & 2.46547 & 2.35487 & 2.46881 & 2.25680 & 2.57155 \\
   \hline
\end{tabular} 
  \caption{Mean values of $T_x, T_y, T_z$ and $m4_x, m4_y, m4_z$ at $t=2$ which are computed using $M_s=100$ standard DSMC simulations with $N=10^8$ particles in each. The initial condition \eqref{eq:IC} is set with the given initial temperatures $[T_x^0 \; T_y^0 \; T_z^0]$.}\label{Table2}
\end{table}

Using the values in~\Cref{Table2}, we can estimate the gradients $\frac{\partial \cJ_1}{\partial \alpha}(\alpha_0)$ via the central finite difference 
\begin{equation}~\label{eq:FD_gradient}
\frac{\partial \cJ_1}{\partial \alpha}(\alpha_0) \approx \frac{\cJ_1(\alpha_0+\Delta \alpha)-\cJ_1(\alpha_0-\Delta \alpha)}{2 \Delta \alpha}.
\end{equation}
We can also estimate their errors $e_{\frac{\partial \cJ_1}{\partial \alpha}(\alpha_0)}$, as shown in~\Cref{Table3}. The error in $\frac{\partial \cJ_1}{\partial \alpha}(\alpha_0)$ consists of two parts:
\begin{equation}~\label{eq:total_error}
e_{\frac{\partial \cJ_1}{\partial \alpha}(\alpha_0)} = e^{FD}_{\frac{\partial \cJ_1}{\partial \alpha}(\alpha_0)} + e^{rand}_{\frac{\partial \cJ_1}{\partial \alpha}(\alpha_0)},
\end{equation}
where the finite difference error can be estimated as 
\begin{equation}~\label{eq:FD_error}
e^{FD}_{\frac{\partial \cJ_1}{\partial \alpha}(\alpha_0) } \approx \frac{\partial^3\cJ_1}{\partial \alpha^3}(\alpha_0) \frac{(\Delta \alpha)^2}{6} \approx 
\frac{-\cJ_1(\alpha_0-2\widetilde{\Delta \alpha})+2\cJ_1(\alpha_0-\widetilde{\Delta \alpha})-2\cJ_1(\alpha_0+\widetilde{\Delta \alpha})+\cJ_1(\alpha_0+2\widetilde{\Delta \alpha})}{12 (\widetilde{\Delta \alpha})^3} (\Delta \alpha)^2,
\end{equation}
where $\widetilde{\Delta \alpha}$ is a finite step used to estimate $\frac{\partial^3\cJ_1}{\partial \alpha^3}$, and the random error can be estimated roughly as
\begin{equation}~\label{eq:rand_error}
e^{rand}_{\frac{\partial \cJ_1}{\partial \alpha}(\alpha_0) } \approx \frac{e^{rand}_{\cJ_1(\alpha_0+\Delta \alpha)}  + e^{rand}_{\cJ_1(\alpha_0-\Delta \alpha)}}{2 \Delta \alpha} \approx \frac{e^{rand}_{\cJ_1(\alpha_0)}}{\Delta \alpha}.
\end{equation}
Notice that the finite difference error is proportional to $(\Delta \alpha)^2$ while the random error is proportional to $1/\Delta \alpha$, meaning that there is an optimal value of 
\begin{equation}~\label{eq:opt_alpha}
\Delta \alpha^* = \left( \frac{3 e^{rand}_{\cJ_1(\alpha_0)} }{ \frac{\partial^3\cJ_1}{\partial \alpha^3}(\alpha_0)} \right)^{1/3} 
\end{equation}
for a given value of the random error $e^{rand}_{\cJ_1(\alpha_0)}$ in $J(\alpha_0)$ and the third-order derivative of $\cJ_1$ at $\alpha_0$ (unknown a priori) for which the total error is minimal. 
We compute the third-order derivatives according to \eqref{eq:FD_error} with $\widetilde{\Delta \alpha}=0.05$. Together with the previously computed $e^{rand}_{T_l}$ and $e^{rand}_{m4_l}$, we estimate the optimal step size $\Delta \alpha^*\in (0.06,1.5)$, for various components of $\frac{\partial T_l}{\partial T_p^0}$ and $\frac{\partial m4_l}{\partial T_p^0}$. Eventually, we opt for $\Delta \alpha=\Delta T_p^0=0.1$ for the computation of the gradients and the estimations of the errors that are shown in~\Cref{Table3}. 
The random errors reflected in the last column of \Cref{Table3} are computed using the standard deviations $\sigma_{T_l}$ and $\sigma_{m4_l}$  from the last column of~\Cref{Table1} together with~\eqref{eq:E[X]approx} and~\eqref{eq:rand_error}. 
Comparing the random errors in the last column with the finite-difference errors in the other columns of~\Cref{Table3}, we observe that our choice of $\Delta \alpha=0.1$ is nearly optimal for $p=x$ components and smaller than optimal for the $p=y,z$ components, but the current value of $\Delta \alpha$ is good enough for our purposes.

\begin{table} 
\small
\setlength{\tabcolsep}{4pt}
\centering 
\begin{tabular}{|c|c|c|c|c|} 
  \hline
\rule{0pt}{10pt}  & $\frac{\partial \cJ_1}{\partial T_x^0} \pm e^{FD}_{\frac{\partial \cJ_1}{\partial T_x^0}}$ & $\frac{\partial \cJ_1}{\partial T_y^0} \pm e^{FD}_{\frac{\partial \cJ_1}{\partial T_y^0}}$ & $\frac{\partial \cJ_1}{\partial T_z^0} \pm e^{FD}_{\frac{\partial \cJ_1}{\partial T_z^0}}$ & $e^{rand}_{\frac{\partial \cJ_1}{\partial T_p^0}}$ \\ \hline
$\cJ_1=T_x(t=2)$  &  0.572335 $\pm$  0.0001  &  0.213839 $\pm$  2.1e-06  &  0.213836 $\pm$  8.6e-07 & $\pm$  0.0002 \\ 
$\cJ_1=T_y(t=2)$  &  0.213831 $\pm$  0.0003  &  0.572331 $\pm$  2.1e-06  &  0.213842 $\pm$  3.1e-07 & $\pm$  0.0002 \\
$\cJ_1=T_z(t=2)$  &  0.213834 $\pm$  0.0002  &  0.213830 $\pm$  2.8e-08  &  0.572322 $\pm$  5.6e-07 & $\pm$  0.0002 \\
$\cJ_1=m4_x(t=2)$ &  2.289942 $\pm$  0.0011  &  0.996551 $\pm$  7.5e-06  &  0.996538 $\pm$  6.0e-06 & $\pm$  0.0009 \\ 
$\cJ_1=m4_y(t=2)$ &  1.071492 $\pm$  0.0025  &  3.147580 $\pm$  7.5e-06  &  1.139506 $\pm$  5.5e-06 & $\pm$  0.0011 \\ 
$\cJ_1=m4_z(t=2)$ &  1.071524 $\pm$  0.0006  &  1.139418 $\pm$  1.1e-06  &  3.147446 $\pm$  5.7e-07 & $\pm$  0.0014 \\ 
   \hline
\end{tabular}
  \caption{Gradients $\frac{\delta \cJ_1}{\delta \alpha}$ that are computed using the finite difference formula \eqref{eq:FD_gradient}, estimates of the corresponding finite-difference errors based on~\eqref{eq:FD_error} and~\Cref{Table2} values, and the corresponding random errors that are computed based on~\eqref{eq:rand_error} and~\Cref{Table1} values. The objective functions $\cJ_1=\{T_l(t=2),m4_l(t=2)$\} and parameter $\alpha=T_p^0$, $l,p\in\{x,y,z\}$. The random errors are the same for all $\alpha=T_p^0$, $p\in\{x,y,z\}$. All quantities here are computed using $\Delta \alpha=0.1$.}\label{Table3}
\end{table}



\subsection{Adjoint DSMC}
In the adjoint DSMC method for the gradient calculation, we first solve the Boltzmann equation~\eqref{eq:homoBoltz} as in the previous examples, but  
at each time step $t=t_k$ the indices of the particles that collided, their collision angles $\sigma_{k,i}$ and the pre-collision relative unit velocities $(v_{k,i} - {\tilde v}_{k,i}){\hat ~}$ are stored in memory for the later use in the backward solve. 

Afterwards, we solve the adjoint DSMC equations derived in~\Cref{DiscreteFormulation} backward in time. Similar to the definition of the forward particles~\eqref{eq:vk}, at the $k$-th time interval, we represent the adjoint particles as
\begin{equation*}
\bGamma_{k} =\{ \bgamma_{1}, \ldots, \bgamma_{i} ,\ldots, \bgamma_{N}\}(t_{k}),
\end{equation*}
and we denote the $i$-th adjoint particle in $\bGamma_k$ as $\bgamma_{k,i}$. We remark that $\bgamma_{k,i}$ is a vector in $\Rthree$. 
We initialize our backward solve with the ``final conditions"~\eqref{gammaFinal}, $\bgamma_{F,i} = - \partial_v r(v_{F,i})$, $i\in 1\dots N$, at the final time $t=T=2$, and $r(v)$ is set to be $v_l^2$ or $v_l^4$, $l\in\{x,y,z\}$, as discussed previously, where $l$ is fixed for a given simulation. Thus, $\bgamma_{F,i}^j = -  2 v_{F,i}^l \delta_{l,j}$ and $\bgamma_{F,i}^j = -  4 (v_{F,i}^l)^3 \delta_{l,j}$, respectively, where $l,j\in\{x,y,z\}$ and $\delta_{l,j}=1$ if $l=j$ and $0$ otherwise. 

The adjoint DSMC equations~\eqref{gammaEqtn}-\eqref{gammaTildeEqtn} are then solved backward in time for each collision that happened during the forward solve using the collision angles and the relative unit velocities that have been stored during the forward solve. Once the backward time evolution reaches the initial time $t=0$, we use \eqref{DalphaEqtn}, $\frac{\partial \cJ_1}{\partial \alpha}=- \frac{1}{N}\sum_{i=1}^N \bgamma_{{I}, i} \cdot \partial_\alpha v_{0,i}(\alpha)$, to calculate the gradient of the objective function with respect to the model parameter $\alpha$. Here, $\partial_\alpha v_{0,i}(\alpha)$ is a derivative of the initial sample of particles with respect to the parameter $\alpha$. Hence, it is important to have not only an initial distribution function $f_0(v,\alpha)$ that depends continuously on the parameter $\alpha$ but also a particular way of sampling such that \textit{the samples depend continuously on the parameter $\alpha$}. Since our initial distribution is an isotropic Gaussian \eqref{eq:IC}, we can sample it by sampling $3N$ values from the standard normal distribution $\cN(0,1)$ and then rescaling the values with appropriate initial temperatures as 
$$v_{0,i}=(v_{0,i}^x,v_{0,i}^y,v_{0,i}^z) = (\sqrt{T_x^0}v_{0,i}^{x\cN},\sqrt{T_y^0}v_{0,i}^{y\cN},\sqrt{T_z^0}v_{0,i}^{z\cN}),$$ 
where $v_{0,i}^{x\cN},v_{0,i}^{y\cN},v_{0,i}^{z\cN}$ are samples of $\cN(0,1)$. For $\alpha=T_p^0$, we can easily compute 
$$\frac{\partial v_{0,i}^j}{\partial T_p^0} = \frac{ v_{0,i}^{j\cN} }{2\sqrt{T_p^0}} \delta_{j,p}=\frac{ v_{0,i}^j}{2T_p^0} \delta_{j,p},\quad j,p\in\{x,y,z\}.$$


\begin{table}
\small
\setlength{\tabcolsep}{4pt}
\centering 
\begin{tabular}{|c|c|c|c|} 
  \hline
\rule{0pt}{10pt}  & $\frac{\partial \cJ_1}{\partial T_x^0} \pm e^{rand}_{\frac{\partial \cJ_1}{\partial T_x^0}}$ & $\frac{\partial \cJ_1}{\partial T_y^0} \pm e^{rand}_{\frac{\partial \cJ_1}{\partial T_y^0}}$ & $\frac{\partial \cJ_1}{\partial T_z^0} \pm e^{rand}_{\frac{\partial \cJ_1}{\partial T_z^0}}$ \\
  \hline
$\cJ_1=T_x(t=2)$  &  0.572316 $\pm$  1.1e-05  &  0.213836 $\pm$  8.4e-06  &  0.213835 $\pm$  7.7e-06  \\
$\cJ_1=T_y(t=2)$  &  0.213846 $\pm$  9.5e-06  &  0.572337 $\pm$  1.1e-05  &  0.213839 $\pm$  9.5e-06  \\
$\cJ_1=T_z(t=2)$  &  0.213838 $\pm$  8.9e-06  &  0.213828 $\pm$  7.5e-06  &  0.572325 $\pm$  1.2e-05  \\
$\cJ_1=m4_x(t=2)$ &  2.289879 $\pm$  1.3e-04  &  0.996589 $\pm$  8.6e-05  &  0.996541 $\pm$  8.3e-05  \\
$\cJ_1=m4_y(t=2)$ &  1.071577 $\pm$  9.1e-05  &  3.147648 $\pm$  1.8e-04  &  1.139492 $\pm$  9.0e-05  \\
$\cJ_1=m4_z(t=2)$ &  1.071486 $\pm$  8.0e-05  &  1.139424 $\pm$  8.0e-05  &  3.147454 $\pm$  1.7e-04  \\
   \hline
\end{tabular} 
  \caption{Gradients $\frac{\delta \cJ_1}{\delta \alpha}$ computed by the adjoint DSMC approach and estimations of the corresponding random errors based on formula~\eqref{eq:E[X]approx}, where the objective function $\cJ_1=\{T_l(t=2),m4_l(t=2)$\} and the parameter $\alpha=T_p^0$, $l,p\in\{x,y,z\}$.}\label{Table4}
\end{table} 

We perform $M_s=100$ forward DSMC and backward adjoint DSMC simulations with $N=10^8$ particles for each simulation. Based on~\eqref{eq:E[X]approx}, we compute the mean values and random errors  for quantities $\frac{\partial T_l}{\partial T_p^0}$ and $\frac{\partial m4_l}{\partial T_p^0}$, $l,p\in\{x,y,z\}$. The results are gathered in~\Cref{Table4}, where the values of the gradients match those in~\Cref{Table3} up to 4-5 significant digits. The random errors in~\Cref{Table4} are approximately $10^{-5}$ for 
$\frac{\partial T_l}{\partial T_p^0}(t=2)$ and $10^{-4}$ for $\frac{\partial m4_l}{\partial T_p^0}(t=2),$
which is about one order of magnitude smaller than the ones in~\Cref{Table3}. It demonstrates another advantage of this method that it does not suffer from amplified random errors compared to the finite difference calculations. For the latter, the random error increases while computing the gradients due to the subtraction of two close values in \eqref{eq:FD_gradient}.



\subsection{DSMC-type scheme for the adjoint equation}~\label{sec:numeric_DSMCtype_method} 
For the DSMC-type scheme, we need to store the same information during the forward modeling as the adjoint DSMC method.
We solve the continuous adjoint equation~\eqref{backwards1} backward in time using a \textit{DSMC-type scheme} that we have developed in~\Cref{sec:DSMCtype_method}. In this method, we evolve values of the function $\gamma(v,t)$ evaluated at locations $v_{k,i}$, i.e., the DSMC particles from the forward solve. At the $k$-th time step, the set
\begin{equation*}
\Gamma_{k} =\{\hat \gamma_1, \ldots , \hat \gamma_i,\ldots \hat \gamma_N\} (t_k)
\end{equation*}
contains all those function values, where $\hat \gamma_i(t_k)=\gamma(v_{k,i},t_k)$ is $i$-th adjoint ``particle" in $\Gamma_k$. We remark that $\hat \gamma_i(t_k)$ here is a scalar value, $\hat \gamma_i(t_k) \in \bbR$. 
We initialize our backward solve with the ``final conditions"~\eqref{backwards0}. That is, $\hat \gamma_i(t=T) = - r(v_{F,i})$, $i\in 1\dots N$, where $\{v_{F,i}\}_{i=1}^N$ are the velocity particles from the forward DSMC at the final time $t=T=2$.
Then the continuous adjoint equation \eqref{backwards1} can be solved backward in time following~\Cref{alg:adj_DSMC}, where $\hat \gamma_i(t_k)$ are updated at $k+1 \rightarrow k$ time step according to \eqref{eq:OTD_DSMC3} only if $v_{k,i}$ particles participated in a collision during the forward solve at the step $k \rightarrow k+1$. To approximate the last terms in~\eqref{eq:OTD_DSMC3}, we perform a linear interpolation using the scattered data $\Gamma_{k+1}$ at locations $v_{k+1,i}$; see~\eqref{eq:OTD_DSMC4}. 
Eventually, once the backward time evolution reaches the initial time $t=0$, the gradient of the objective function~\eqref{optimization1} can be computed using the values of $\hat \gamma_i(t=0)$, i.e., $\Gamma_0$.
In case of the initial condition \eqref{eq:IC} and the parameter $\alpha=T_p^0$, $p\in\{x,y,z\}$, the derivative $\partial_\alpha f_0(v;\alpha)$ of the initial distribution function with respect to the parameter $\alpha$ in~\eqref{optimization1} becomes  
\begin{equation} \label{eq:IC_derivative}
\frac{\partial f_0(v;\alpha)}{\partial \alpha}=\left(\frac{v_p^2}{T_p^0} - 1\right)\frac{1}{2T_p^0} f_0(v).
\end{equation}
By approximating the initial distribution function with the sample particles $f_0(v) \approx  \frac{1}{N} \sum_{i=1}^N \delta(v-v_{0,i})$, we get
\begin{equation} \label{eq:gradient_DSMCtype}
\frac{\partial J_1}{\partial \alpha} = -\int \gamma(v,0)  \frac{\partial f_0(v;\alpha)}{\partial \alpha} dv =-\int \gamma(v,0) \left(\frac{v_p^2}{T_p^0} - 1\right)\frac{1}{2T_p^0} f_0(v)dv \approx -\frac{1}{N}  \sum_{i=1}^N \frac{\hat \gamma_i(t=0)}{2T_p^0}  \left(\frac{(v_{0,i}^p)^2}{T_p^0} - 1\right)  = \frac{\partial \cJ_1}{\partial \alpha}.
\end{equation}

\begin{table} 
\small
\setlength{\tabcolsep}{4pt}
\centering 
\begin{tabular}{|c|c|c|c|} 
  \hline
\rule{0pt}{10pt}  & $\frac{\partial \cJ_1}{\partial T_x^0} \pm e^{rand}_{\frac{\partial \cJ_1}{\partial T_x^0}}$ & $\frac{\partial \cJ_1}{\partial T_y^0} \pm e^{rand}_{\frac{\partial \cJ_1}{\partial T_y^0}}$ & $\frac{\partial \cJ_1}{\partial T_z^0} \pm e^{rand}_{\frac{\partial \cJ_1}{\partial T_z^0}}$ \\
  \hline
$\cJ_1=T_x(t=2)$  &  0.572229 $\pm$  3.6e-03  &  0.213398 $\pm$  8.9e-04  &  0.213809 $\pm$  1.1e-03  \\
$\cJ_1=T_y(t=2)$  &  0.211871 $\pm$  4.3e-03  &  0.573324 $\pm$  2.6e-03  &  0.213879 $\pm$  1.7e-03  \\
$\cJ_1=T_z(t=2)$  &  0.216833 $\pm$  3.2e-03  &  0.211929 $\pm$  1.3e-03  &  0.572561 $\pm$  2.3e-03  \\
$\cJ_1=m4_x(t=2)$ &  2.285860 $\pm$  1.9e-02  &  0.992066 $\pm$  8.6e-03  &  0.992431 $\pm$  9.4e-03  \\
$\cJ_1=m4_y(t=2)$ &  1.066416 $\pm$  2.7e-02  &  3.147584 $\pm$  1.7e-02  &  1.133434 $\pm$  1.4e-02  \\
$\cJ_1=m4_z(t=2)$ &  1.079131 $\pm$  1.9e-02  &  1.127308 $\pm$  9.9e-03  &  3.145782 $\pm$  2.7e-02  \\
   \hline
\end{tabular} 
  \caption{Gradients $\frac{\delta \cJ_1}{\delta \alpha}$ computed by the DSMC-type scheme for the adjoint equation and estimations of the corresponding random errors based on formula~\eqref{eq:E[X]approx}, where the objective function $\cJ_1=\{T_l(t=2),m4_l(t=2)\}$ and the parameter $\alpha=T_p^0$, $l,p\in\{x,y,z\}$.}\label{Table5}
\end{table}

We performed $M_s=10$ forward DSMC simulations and backward simulations of the adjoint equation using the DSMC-type scheme with $N=10^6$ particles for each simulation. We choose smaller $N$ and $M_s$ here because this scheme works much slower due to the interpolation procedure at every time step for every colliding particle. We compute the mean values of the gradients according to \eqref{eq:gradient_DSMCtype}, and estimate their random errors based on~\eqref{eq:E[X]approx}. 
The results are presented in~\Cref{Table5}. We observe that the gradients in \Cref{Table4} and \Cref{Table5} match up to 2-4 significant digits. The random errors in~\Cref{Table5} are approximately 0.003 for $\frac{\partial T_l}{\partial T_p^0}(t=2)$ and 0.02 for $\frac{\partial m4_l}{\partial T_p^0}(t=2),$ which are about two orders of magnitude larger than the ones in~\Cref{Table4} as a result of using $N=10^6$ instead of $N=10^8$ and the additional error introduced by numerical interpolation.



\subsection{Direct discretization of the adjoint equation integrals} \label{sec:Direct_num}
Here we solve the continuous equation~\eqref{backwards1} for $\gamma(v,t)$ from $t=T=2$ to $t=0$ with the ``final condition" \eqref{backwards0} using a direct numerical integration scheme to compute the integral term on the right-hand side. See details of the numerical scheme in~\ref{A1}.

As previously, we consider the gradients of the second-order moments, $T_l(t=2)$, and the fourth-order moments, $m4_l(t=2)$, with respect to $\alpha=T_p^0$, $l,p\in\{x,y,z\}$.
For each objective function, we perform one forward DSMC simulation with $N=10^8$ particles and one backward simulation via the scheme described in~\eqref{backwards1_simplified2_discrete}. We did only one simulation ($M_s=1$) for each objective function due to its extensively long computational time; see the next subsection for discussions on performance. The results are gathered in~\Cref{Table6}. Values of the gradients in~\Cref{Table6} and~\Cref{Table4} match up to 2-3 significant digits. 
We still have random errors that are contributed through the values of $f(v,t)$ in the forward DSMC. The random errors in~\Cref{Table6} can be roughly estimated using two standard deviations in the last column of \Cref{Table1}, or $e^{rand}_{T_l} = 2 \sigma_{T_l}/\sqrt{M_s} \approx 0.0002$ and $e^{rand}_{m4_l} = 2 \sigma_{m4_l}/\sqrt{M_s} \approx 0.001$. They are much smaller than the overall errors (dominated by finite discretization errors here) 0.002 for $\frac{\partial T_l(t=2)}{\partial T_p^0}$  and 0.05 for $\frac{\partial m4_l(t=2)}{\partial T_p^0},$
which we can compute by comparing the gradient values in ~\Cref{Table6} and~\Cref{Table4}.

We have also numerically verified Equation~\eqref{eq:dJ4} using the values of $\gamma(v,t_k)$ at grid points $v=v_{i_x,i_y,i_z}$ obtained during the direct integration simulations and using the finite difference to compute the derivative $\gamma'(v,t_k)$. We approximate $\mathbb E[\bgamma_{k,i}|v=v_{k,i}]$ using the histogram count of $\bgamma_{k,i}$ obtained via the adjoint DSMC method around the same grid points $v=v_{i_x,i_y,i_z}$.

\begin{table} 
\small
\centering 
\begin{tabular}{|c|c|c|c|} 
  \hline
   & $\frac{\partial J_1}{\partial T_x^0}$ & $\frac{\partial J_1}{\partial T_y^0}$ & $\frac{\partial J_1}{\partial T_z^0}$ \\
  \hline
$J_1=T_x(t=2)$  &  0.573407   &  0.214738   &  0.214742  \\ 
$J_1=T_y(t=2)$  &  0.214944   &  0.572962   &  0.214757  \\
$J_1=T_z(t=2)$  &  0.211709   &  0.211531   &  0.569767  \\
$J_1=m4_x(t=2)$ &  2.343578   &  1.027094   &  1.026818  \\
$J_1=m4_y(t=2)$ &  1.099767   &  3.193492   &  1.165838  \\
$J_1=m4_z(t=2)$ &  1.072768   &  1.137018   &  3.153235  \\
   \hline
\end{tabular} 
\caption{Gradients $\frac{\delta J_1}{\delta \alpha}$ computed by the direct discretization of the integrals in the continuous adjoint equation~\eqref{backwards1}, where $J_1=\{T_l(t=2),m4_l(t=2)$\}, $\alpha=T_p^0$, $l,p\in\{x,y,z\}$.}\label{Table6}
\end{table}



\subsection{Method comparison}
So far, we have described and demonstrated four different ways to compute the gradient of an objective function numerically after the forward DSMC simulations. The summary of their comparison in terms of memory requirements, error scaling and operation count is given in~\Cref{Table7}. See the details of the comparison in~\ref{A2}.

\Cref{Table8} shows the timings of our code we recorded per objective function using Intel Core i7-3770K processor (4 cores @4.5Ghz) and $N=10^6,10^7,10^8$, $T=2$, $\Delta t=0.1$, $n_\text{grid}=30,40, n_{\varphi}=n_{\theta}=10$ parameters.

\begin{table} 
\footnotesize
\setlength{\tabcolsep}{4pt}
\centering 
\begin{tabular}{|m{3.2cm}|m{3.3cm}|m{3.0cm}|m{3.1cm}|} 
  \hline
   &  memory requirements, bytes  & error scaling & operation count \\
  \hline 
(i) finite difference                         &    $24N$                                  & $\cO((\Delta \alpha)^2)+\cO(\frac{1}{\sqrt{N} \Delta \alpha})$    &  $23.5 (D_\alpha+1) N\mu T$  \\
(ii) adjoint DSMC method                                      &  $24N$ + $(24N + 28 N\mu T)$                    & $\cO(\frac{1}{\sqrt{N}})$              & $17.5 N\mu T$  \\
(iii) DSMC-type scheme & $24N$ + $(8N + 28 N\mu T)$                                   & $\cO(\frac{1}{\sqrt{N}})$              & $\approx (40 + \frac{C}{\Delta t} \log(N)) N\mu T$  \\
(iv) direct discretization &  $8 n_\text{grid}^3 (\frac{T}{\Delta t}+2)$  & $\cO(\frac{1}{\sqrt{N}})$ + $\cO((\Delta v)^2)$  & $\approx (21 n_{\varphi}n_{\theta}+4)n_\text{grid}^6 \frac{T}{\Delta t}$  \\
   \hline
\end{tabular}
  \caption{The comparison of the four methods that compute the gradients of the Boltzmann-constrained optimization problems, in terms of memory requirements, error scaling and operation count.}  \label{Table7}
\end{table}

\begin{table} 
\centering
{\renewcommand{\arraystretch}{1.2}%
\begin{tabular}{|c|c|c|c|c|c|}   \hline 
 \quad  &  $N = 10^6$  & $N = 10^7$ & $N = 10^8$ & $n_\text{grid}=30$& $n_\text{grid}=40$ \\
  \hline 
forward DSMC simulation (Algorithm~\ref{alg:DSMC}) &    $0.38$   sec & $5$ sec &  $60$ sec  & &\\
adjoint DSMC simulation (Algorithm~\ref{alg:adj_DSMC1}) &  $0.22$ sec & $2.7$ sec & $30$ sec & & \\
the DSMC-type scheme (Algorithm~\ref{alg:adj_DSMC}) &  $280$  sec  & $4100$ sec &\quad & & \\
direct discretization of the equation~\eqref{backwards1} &   &  & & $25000$ sec & $126000$ sec\\
   \hline
\end{tabular} }
  \caption{CPU run time of the four different methods presented in Section~\ref{sec:Forward_DSMC}-\ref{sec:Direct_num} under different parameters.} \label{Table8}
\end{table}

The adjoint DSMC is slightly faster (up to 20-30$\%$) than the forward DSMC for the same number of particles $N$ (timings for computing the vector $\sigma$ and the unit collision direction are included in the forward DSMC timings above and take about $20\%$ of the overall timings, whereas they would add about $35\%$ to the adjoint DSMC timings if we were to include them there). The adjoint DSMC is more than 1000 times faster than the DSMC-like scheme 
and much faster than the direct discretization of~\eqref{backwards1}. At the same time, the errors in the adjoint DSMC are at least one order of magnitude smaller than in other methods due to absence of finite-difference or interpolation errors; see Tables \ref{Table3}-\ref{Table5} and \Cref{sec:Direct_num} where we have error estimates (partly numeric and partly analytic). By performing only one adjoint DSMC simulation, we obtain all the gradient components of the objective function. At the same time, if we use the finite difference method, we need to perform at least $D_\alpha+1$ simulations, where $D_\alpha$ is the dimensionality of $\alpha$. The benefits of the adjoint DSMC algorithm particularly stand out when solving large-scale optimization problems.

\begin{figure}
\centering
\includegraphics[width=0.6\textwidth]{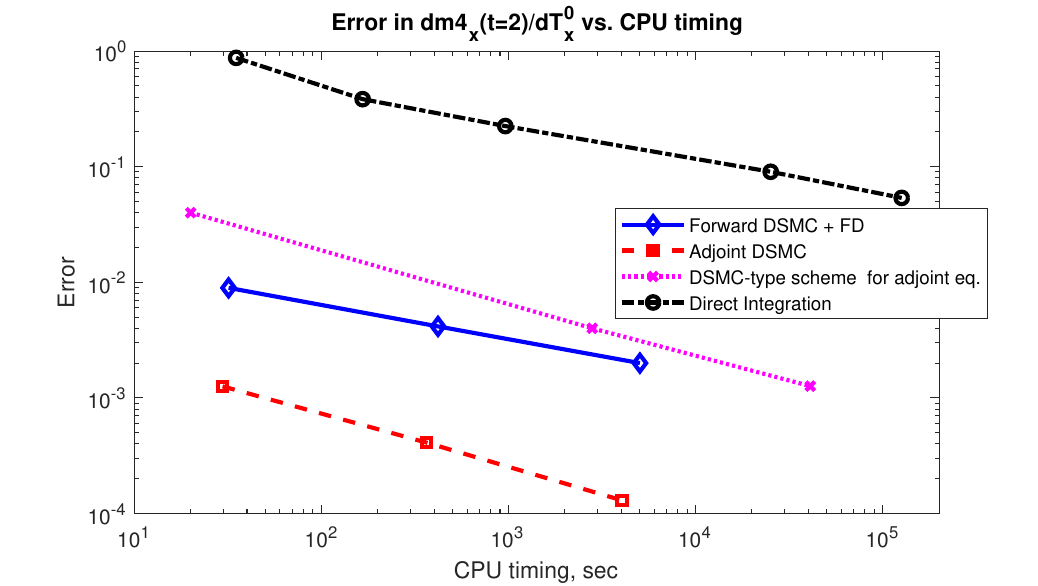}
\caption{The error in $\frac{\partial m4_x(t=2)}{\partial T_x^0}$ vs. CPU time for all four methods} \label{fig:Error_vs_CPU}
\end{figure}

Fig.~\ref{fig:Error_vs_CPU} shows the error in $\frac{\partial m4_x}{\partial T_x^0}$ vs.~CPU time measured in simulations for all four methods. Since one forward DSMC solve is needed for all methods, Fig.~\ref{fig:Error_vs_CPU} reflects only CPU timing for the additional computations needed to compute the gradient, namely one extra forward DSMC simulation ($D_\alpha=1$) for (i) or one backward solution for (ii)-(iv). 
For (ii), the errors were estimated using formula~\eqref{eq:E[X]approx} as the error values in \Cref{Table4}.
For (iii) and (iv) the error values in the figure were computed as a difference between the numerical values $\frac{\partial m4_x}{\partial T_x^0}$ obtained in numerical simulations and the reference value of $\frac{\partial m4_x}{\partial T_x^0}$ from \Cref{Table4} as it is the most accurate value we have computed. 
For (i) the errors were estimated using formulas \eqref{eq:total_error}-\eqref{eq:opt_alpha} and the fact that $e^{rand}_{m4_x(t=2)} \propto 1/\sqrt{N}$ to demonstrate that the error scaling and the actual errors (in comparison to the reference value from \Cref{Table4}) were $2-10$ times smaller.
Note that the slope of the forward DSMC+FD line (blue) is less than the slopes for the adjoint DSMC (red) and DSMC-type scheme for adjoint equation (purple) since the total error in~\eqref{eq:total_error} for the optimal $\Delta\alpha^*$ (as in \eqref{eq:opt_alpha}) scales like $\propto 1/N^{1/3}$ and CPU time $\propto N$.




\subsection{Optimization examples}
We have previously discussed the accuracy and performance of several different numerical schemes to compute the gradient of optimization problems constrained by the Boltzmann equation. The adjoint DSMC method particularly stands out for its simplicity, computational efficiency, and the direct connections with the well-established forward DSMC method~\cite{nanbu1980direct}, as discussed in~\Cref{sec:DSMC_property}. Here, we use two optimization examples to illustrate the great potential of the adjoint DSMC method for efficiently solving optimization problems constrained by the Boltzmann equation with the nonlinear collision operator.

\subsubsection{Matching the velocity moments}\label{sec:opt_ex1}
We have been using the velocity moments of the probability distribution at the final time $T$ as the objective function to test the accuracy of the gradients. Here, we follow the earlier discussions and set the first objective function as
\begin{equation}\label{eq:opt_ex1_obj}
\cJ_1(\alpha)= \|\bd_1 - \bd_2 \|_2^2
\end{equation}
where $\bd_1 = [T_x, T_y, T_z]$, the second velocity moments in each direction and $\bd_2 = \frac{1}{2}[m4_x, m4_y, m4_z]$, half of the fourth-order velocity moments in each direction at $t=T=2$. As we have defined earlier, $$m4_l(T)=\frac{1}{N}\sum_{i=1}^N{(v_{F,i}^l)^4 } \approx \int_{\bbR^3} v_l^4 f(v,T) dv,\  T_l(T)=\frac{1}{N}\sum_{i=1}^N{(v_{F,i}^l)^2 } \approx \int_{\bbR^3} v_l^2 f(v,T) dv,$$
for $l\in \{x,y,z\}$. Here, $f(v,T)$ solves the Boltzmann equation~\eqref{eq:homoBoltz} given the initial condition~\eqref{eq:IC}. We only treat the initial temperature of the $y$ direction, $T_y^0$, as the unknown parameter $\alpha$, while fixing $T_x^0 = 0.5$ and $T_z^0 = 1$. This is to avoid the trivial optimal solution $T_x^0 = T_y^0 = T_z^0 = 0$ that minimizes the objective function~\eqref{eq:opt_ex1_obj} if $\alpha = [T^0_x, T^0_y, T^0_z]$.

One may notice that the objective function~\eqref{eq:opt_ex1_obj} in this example does not match the formulation~\eqref{eq:OTD_obj}. Nevertheless, it is easy to adapt a general objective function to either the continuous or the DSMC adjoint system based on the role of the adjoint equations: large systems of chain rule which propagate the Fr\'echet derivative $\frac{\delta J_1}{\delta f(v,T)}$ backward in time to $\frac{\delta J_1}{\delta f(v,0)}$ as $f(v,0)$ directly depends on the model parameter $\alpha$. Hence, in this example, the final condition for the continuous adjoint equation~\eqref{backwards1} should be
$$
\gamma(v,T) = - \frac{\delta J_1}{\delta f(v,T)} = -2(\bd_1 - \bd_2)\cdot \frac{\delta (\bd_1 - \bd_2)}{\delta f(v,T)}=  -2(\bd_1 - \bd_2)\cdot (v^2-\frac{v^4}{2}).
$$
The final condition for the DSMC adjoint system, $\bgamma_{F,i} = [\bgamma^x_{F,i},\bgamma^y_{F,i},\bgamma^z_{F,i}]$, also the adjoint variable for the final particle velocity $v_{F,i} = [v^x_{F,i} ,v^y_{F,i} ,v^z_{F,i} ]$, should be
$$
\bgamma^l_{F,i} = - \frac{\partial \cJ_1}{\partial v^l_{F,i}} = -\frac{2}{N}(T_l(T) - m4_l(T))(2v^l_{F,i} - 2(v^l_{F,i})^3),\,l\in{x,y,z},\,i = 1,2,\dots,N.
$$

\begin{figure}
\centering
\includegraphics[width=0.49\textwidth]{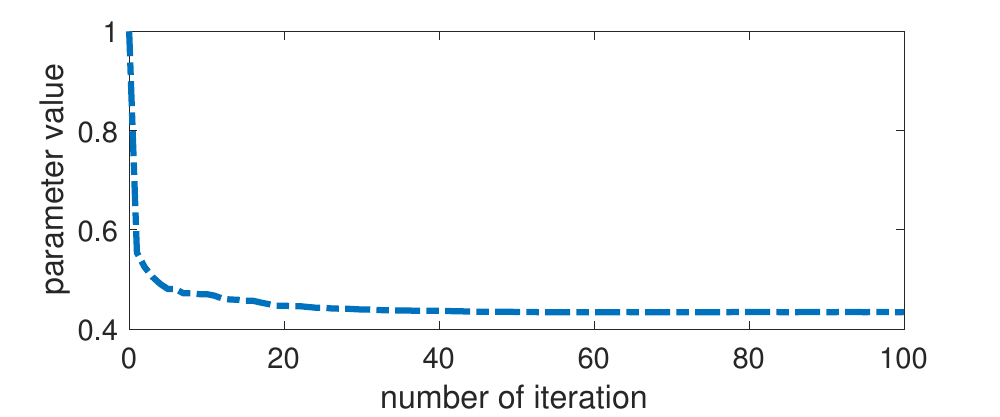}
\includegraphics[width=0.49\textwidth]{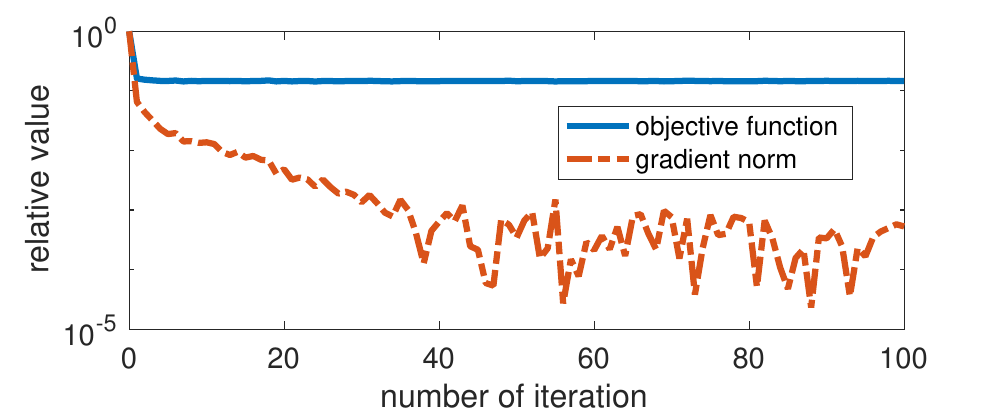}
\caption{Left: the convergence history of the parameter $T_y^0$ for the example discussed in~\Cref{sec:opt_ex1}. Right: the decrease of the normalized objective function value and the size of the gradient for the first 100 iterations in the optimization step.} \label{fig:Optimization_ex1}
\end{figure}

Starting with $T_y^0 = 1$ as the initial guess for $\alpha$, we use a gradient-based optimization algorithm to minimize the objective function~\eqref{eq:opt_ex1_obj}. The steepest descent method with a backtracking line search following the Armijo--Goldstein condition is applied to find a proper stepsize along the descent direction~\cite{nocedal2006numerical}. We compute the gradient by solving one forward DSMC with the current $\alpha$, and then one adjoint DSMC is solved in every iteration of the optimization process. The total number of particles in both DSMC simulations is $N=10^7$. The spacing in the time domain is $\Delta t= 0.1$. The convergence history of this example is shown in Fig.~\ref{fig:Optimization_ex1}. Both the objective function and the size of the gradient monotonically decrease in the first 30 iterations. The convergence slows down as the gradient is smaller than $0.1\%$ of its initial size. We start to observe oscillations in the gradient that come from the random errors in the DSMC solutions. The iterates converge to $0.4344$, the global minimum of the objective function if $T_y^0$ is the only parameter.

\begin{figure}
\centering
\includegraphics[width=0.49\textwidth]{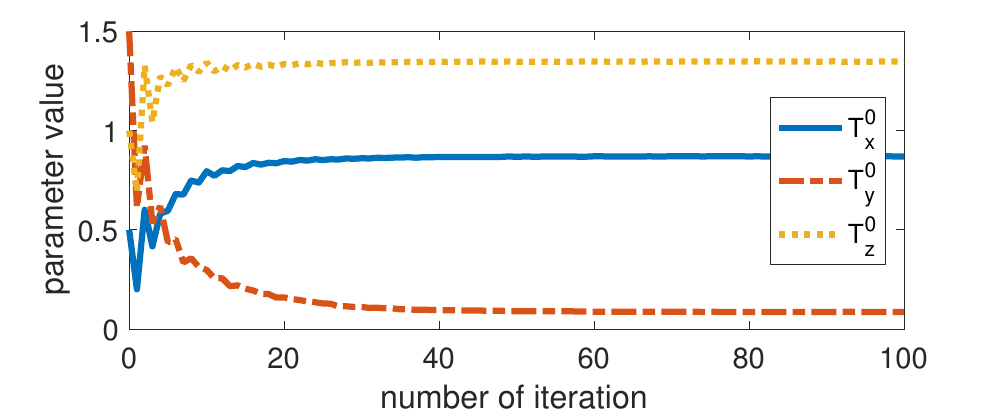}
\includegraphics[width=0.49\textwidth]{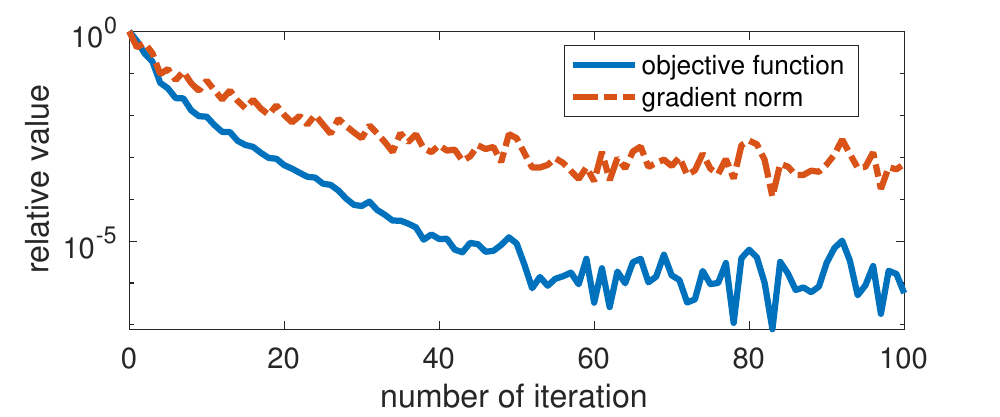}
\caption{Left: the convergence history of the three parameters $[T_x^0,T_y^0,T_z^0]$ for the example discussed in~\Cref{sec:opt_ex2}. Right: the decrease of the normalized objective function value and the $\ell^2$ norm of the gradient for the first 100 iterations in the optimization step.} \label{fig:Optimization_ex2}
\end{figure}

\subsubsection{Inverse Problem}\label{sec:opt_ex2}
Our second example is based the setup of an inverse problem. The fourth-order velocity moments at the final time $T = 2$ are statistical quantities of interest that can be observed in a realistic or experimental setting. The observable information is solely affected by the unknown initial temperature of the distribution that we aim to recover by minimizing the difference between the observed data and the predicted data simulated by our guess of the model parameter. The reconstruction is formulated as a nonlinear least-squares problem
\begin{equation*}
\alpha^* = \argmin\limits_{\alpha} J_1(\alpha)= \argmin\limits_{\alpha}  \|\bd_{\text{obs}} - \bd_\text{pred}(\alpha)\|^2_2
\end{equation*}
where the predicted data $\bd_{\text{pred}}(\alpha) = [m4_x, m4_y ,m4_z]$ is a vector of the fourth-order velocity moments at $T = 2$. We set the true data $\bd_{\text{obs}} = [2,1,3]$. All other notations and the choice of optimization algorithm follow the previous optimization example. 

The initial guess of the parameters is $[0.5,1.5,1.0]$. The convergence history of the computational inverse problem is shown in~Fig.~\ref{fig:Optimization_ex2}. The three components of $\alpha$ converge to $[0.8670, 0.0870, 1.3470]$ in the first $40$ iterations. The same as in the previous example, the small variations in the objective function and the gradient norm for the remaining $60$ iterations are introduced by the random errors of DSMC simulations, as seen in the plots. 
Increasing the number of particles can help mitigate the small perturbations.

\begin{remark}
We show two simple optimization experiments as examples, but one can apply the framework to more general and large-scale optimization problems constrained by the nonlinear Boltzmann equation. There are at least three directions to generalize the applications. First, the dimensionality of the unknown parameter could be increased with hardly any extra cost. Second, the model parameter is not limited to the initial condition. Such examples include shape optimization of the flow channel~\cite{sato2019topology}. Third, following the same idea, we plan to generalize the adjoint DSMC systems for inhomogeneous Boltzmann equation with the nonlinear collision operator or even more complicated kinetic description.
\end{remark}


\section{Conclusion}\label{sec:conclusion}
In this paper, we present the OTD and DTO approaches of computing the gradient of optimization problems based on the nonlinear Boltzmann equation. The highlight of both frameworks is that one only needs to solve the Boltzmann equation and the adjoint system once to compute the gradient, independent of the size of the unknown in the optimization. The adjoint DSMC system, derived by the DTO approach, offers a deterministic numerical scheme that is remarkably efficient to implement with the forward DSMC method. On the other hand, the Monte Carlo type method designed for the continuous adjoint equation could potentially be used for the linear Boltzmann equation~\cite{bobylev2008some}. We plan to extend both frameworks to the general VHS kernel and to the inhomogeneous case. In particular, the adjoint DSMC approach applies to more general kinetic models for gases and plasmas whose behavior could be modeled by Monte Carlo binary collisions. One of such models is the Coulomb collision for charged particles~\cite{wang2008particle}. As the next step, we will apply the adjoint DSMC methods to realistic optimization problems that occur naturally in a broader class of kinetic applications. 

\section*{Acknowledgments}
This material is based upon work supported by the National Science Foundation under Award Number DMS-1913129 and the U.S. Department of Energy under Award Number DE-FG02-86ER53223. The authors thank the Courant Institute of Mathematical Sciences, New York University, for research support and computational resources.

\bibliographystyle{elsarticle-harv}
\bibliography{main}

\appendix

\section{Direct numerical integration scheme of the continuous adjoint equation \eqref{backwards0}}~\label{A1}

We treat $\gamma(v,t)$ as a continuous function with scalar values. We consider a grid in the $v\in\Rthree$ space with $n_\text{grid}$ number of grid points equally spaced in the interval $[-5v_{th}, 5v_{th}]$ for each of the $x,y,z$ directions. We set the thermal velocity $v_{th}=\sqrt{T_M}$, where $T_M=(T_x^0+T_y^0+T_z^0)/3$ is the equilibrium temperature. Thus, the grid spacing is $\Delta v = 10v_{th}/(n_\text{grid}-1)$ and the grid points are $v_{i_x,i_y,i_z}=[-5v_{th}+i_x \Delta v, -5v_{th}+i_y \Delta v, -5v_{th}+i_z \Delta v]$, where $i_x,i_y,i_z \in \{0,1,...,n_\text{grid}-1\}$.

To propagate equation \eqref{backwards1} backward in time, we need values of $f(v,t)$ at the grid points at each time step. Hence, at each time step of the forward DSMC solve of the Boltzmann equation~\eqref{eq:homoBoltz}, we compute and store those function values using a 3D histogram with bins of size $\Delta v$ in each direction that are centered at the grid points $v_{i_x,i_y,i_z}$.

To solve equation~\eqref{backwards1} numerically we first simplify it as follows
\begin{equation} \label{backwards1_simplified}  
-\frac{\partial \gamma(v,t)}{\partial t} = \int_\Rthree \int_\Stwo(\gamma(v_1') + \gamma(v')) f(v_1) q d\sigma dv_1 - \frac{\mu}{\rho} \int_\Rthree \gamma(v_1) f(v_1) dv_1 -  \mu\gamma(v)
\end{equation}
where $v',v_1'$ are as in~\eqref{BoltzmannSolution} and $\sigma$ is a unit vector spanning the unit sphere. Due to the symmetry 
$v'(-\sigma)=v_1'(\sigma)$ and the fact that $\sigma$ spans the whole unit sphere, $\gamma(v_1')$ and $\gamma(v')$ in the first integral give equal contributions.
 

We rewrite the integral over $\sigma$ as an integral over two angles, $\varphi\in[-\pi,\pi]$ and $\theta\in[0,\pi]$, and thus replace $d\sigma$ with $\sin\theta d\theta d\varphi = d(-\cos\theta)d\varphi$.
We discretize $\varphi$ with $n_{\varphi}$ points, and then $\varphi_j = -\pi + 2j\pi/n_{\varphi}$, $j=\{0,1,...,n_{\varphi}-1\}$. Similarly, we discretize $\cos\theta$ with $n_{\theta}$ points, and then $\left(\cos\theta\right)_h = -1 + 1/n_{\theta} + 2h/n_{\theta}$, $h=\{0,1,...,n_{\theta}-1\}$. Thus, 
\[\sigma_{j,h}=(\cos\varphi_j(\sin\theta)_h,\sin\varphi_j(\sin\theta)_h,(\cos\theta)_h), \text{\;where\;} (\sin\theta)_h=\sqrt{1-(\cos\theta)_h^2} \text{\;and\;} \Delta \sigma= 4\pi/(n_{\varphi}n_{\theta}).\]

\begin{figure}
\centering
\includegraphics[width=0.45\textwidth]{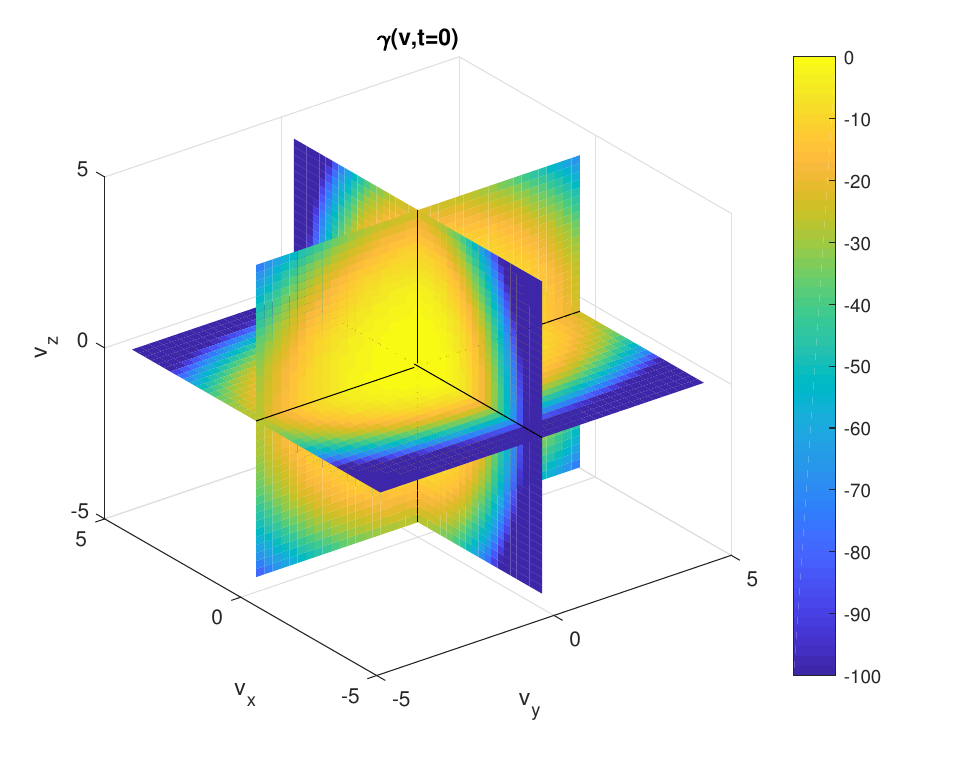}
\includegraphics[width=0.45\textwidth]{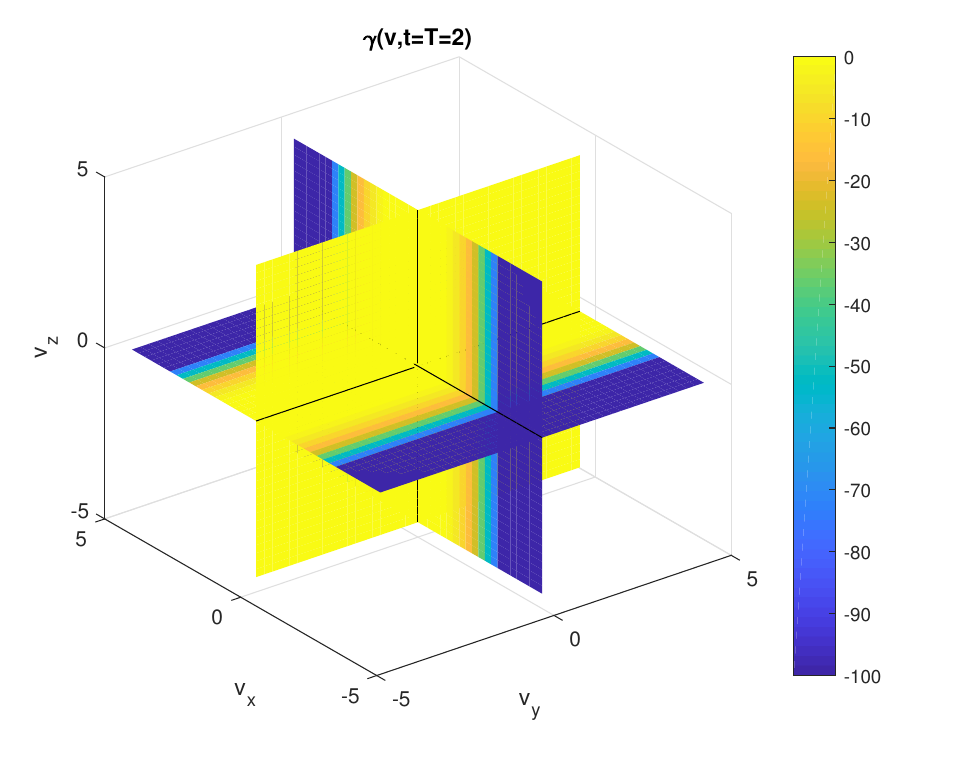}
\caption{Slices of the function $\gamma(v,t)$ at $t=0$ (left) and $t=T=2$ (right) from the backward evolution of equation~\eqref{backwards1_simplified} using the numerical scheme \eqref{backwards1_simplified2_discrete}. Here, $\gamma(v,t=T)=-r(v) = -v_x^4$, $n_\text{grid}=40$, $n_{\varphi}=n_{\theta}=10$ and $\Delta t=0.1$.} \label{fig:gamma_Integrals}
\end{figure}

Finally, we discretize the integrals as sums over the grid $\{v_{i_x,i_y,i_z}, \varphi_j, (\cos\theta)_h\}$, $i_x,i_y,i_z = \{0,1,...,n_\text{grid}-1\}$, $j=\{0,1,...,n_{\varphi}-1\}$, $h=\{0,1,...,n_{\theta}-1\}$ and use the backward Euler scheme for the time derivatives to obtain the following numerical scheme:
\begin{align}
&\quad\frac{\gamma(v_{i_x,i_y,i_z},t_k) - \gamma(v_{i_x,i_y,i_z},t_{k+1}) }{ \Delta t} \label{backwards1_simplified2_discrete} \\ &=2\sum_{i'_x,i'_y,i'_z=0}^{n_\text{grid}-1} \left( \sum_{j=0}^{n_{\varphi}-1} \sum_{h=0}^{n_{\theta}-1} \gamma(v',t_{k+1}) \right) f(v_{i'_x,i'_y,i'_z},t_{k+1}) \frac{\Delta \sigma (\Delta v)^3}{4\pi} - \sum_{i'_x,i'_y,i'_z=0}^{n_\text{grid}-1} \gamma(v_{i'_x,i'_y,i'_z},t_{k+1}) f(v_{i'_x,i'_y,i'_z},t_{k+1}) (\Delta v)^3 -  \gamma(v_{i_x,i_y,i_z},t_{k+1})\nonumber,
\end{align} 
where $v'= 1/2(v_{i_x,i_y,i_z}+v_{i'_x,i'_y,i'_z}) + 1/2|v_{i_x,i_y,i_z}-v_{i'_x,i'_y,i'_z}|\sigma_{j,h}$ and we used our assumptions that $q(\sigma)=1/(4\pi)$ and $\rho =1$, and thus $\mu = \rho \int_\Stwo q(\sigma) d\sigma=1$. Since post-collision velocity $v'$ does not always fall onto the grid points of the velocity space, we approximate the function value $\gamma(v',t_{k+1})$ by linear interpolation using the set of function values $\gamma(v_{i_x,i_y,i_z},t_{k+1})$.
Fig.~\ref{fig:gamma_Integrals} shows the function $\gamma(v,t)$ computed according to the numerical scheme \eqref{backwards1_simplified2_discrete} at $t=T=2$ and $t=0$ using $\gamma(v,t=T)=-r(v) = -v_x^4$ and $n_\text{grid}=40, n_{\varphi}=n_{\theta}=10$, $\Delta t=0.1$.

To compute the gradients of the objective function with respect to the parameter $\alpha$, we use the values of $\gamma(v,0)$ and the derivative of the initial distribution with respect to $\alpha$; see~\eqref{optimization1}. For the initial distribution~\eqref{eq:IC} where the parameter $\alpha=T_p^0$, $p\in\{x,y,z\}$, the derivative of the initial distribution is computed following~\eqref{eq:IC_derivative}.
After approximating the integral in~\eqref{optimization1} by the Riemann sum over the grids, we have
\begin{equation*}
\frac{\partial J_1}{\partial T_p^0} = -\int \gamma(v,0)  \frac{\partial f_0(v;T_p^0)}{\partial T_p^0} dv =-\int \gamma(v,0) \left(\frac{v_p^2}{T_p^0} - 1\right)\frac{1}{2T_p^0} f_0(v) dv  \approx  \sum_{i_x,i_y,i_z=0}^{n_\text{grid}-1}  \gamma(v_{i_x,i_y,i_z},0) \left(\frac{v_{i_p}^2}{T_p^0} - 1\right)\frac{1}{2T_p^0} f_0(v_{i_x,i_y,i_z}) (\Delta v)^3.
\end{equation*}

\section{Numerical methods comparison. Memory requirements, error scaling and operation count.} \label{A2}

\subsection{Memory requirements}
(i) The forward DSMC simulations are done using the Nanbu--Babobsky scheme with $N$ particles. Thus, we need $3N \times 8 = 24 N$ bytes in double-precision arithmetic to store the velocity of the particles. At each time step, there are $N_c= N \mu \Delta t$ particles that collide, which is a small fraction of $N$. Thus, we disregard the temporary variables used for the computation of the collisions in calculating the memory requirements. We also override the particle velocities to not use extra memory for the new velocities at each time step. To compute the gradients by the finite difference, we run several simulations with slightly different values of the parameter $\alpha$. Since these simulations can be done sequentially, it does not increase the overall memory requirements.

(ii) The adjoint DSMC method requires $3N \times 8 = 24N$ bytes for the storage of $\bgamma_{k,i}$ in addition to the memory requirements from the forward DSMC propagation. During the forward DSMC simulation, we also need to store the following information about the collisions that happen at each time step. For each colliding pair, we store
\begin{itemize}
  \item the indices of the two colliding particles, 4 bytes each if stored in the UInt32 format allowing for values in $[0,2^{32}-1]$; since $2^{32} \approx 4.3 \times 10^9 > 10^8$, it is good enough for our purposes;
  \item the $\sigma$ vector in $\Rthree$, 8 bytes for each component of the vector;
  \item the unit collision direction $\frac{v-v_1}{|v-v_1|} \in \Rthree$, 8 bytes for each component.
\end{itemize}

In total, the memory requirements are 56 bytes per colliding pair or 28 bytes per colliding particle. Overall, we need $28N_c$ new bytes stored per time step. If the number of time steps is $M =T/\Delta t$, the total amount of extra storage needed for the backward propagation is $28 N_c M  = 28 N\mu \Delta t  T/\Delta t =  28 N\mu T$. 
It becomes comparable to the $24N$ bytes needed for the particle storage in the forward DSMC when  $28 N\mu T  > 24 N$, or when roughly $\mu T > 1$. The total amount of extra memory required for the backward propagation is $24N + 28 N\mu T$ bytes. In our simulations, we use $\mu=1$, $T=2$, so the total amount of memory required for the adjoint DSMC was about four times the amount of memory required for the forward DSMC.

(iii) The DSMC-type scheme for the continuous adjoint equation~\eqref{backwards1} requires $1N \times 8 = 8 N$ extra bytes for the storage of $\gamma(v_{k,i},t_k)$. Like the forward DSMC, we do not consider temporary storage needed for the colliding particles at each time step since they constitute only a small fraction of the overall number of particles and the new data for $\gamma(v_i,t_k)$ overrides the old one. We also do not store $v_{k,i}$ at every time step during the forward solve since we can always recover the velocity particles when marching back from $t=T$ to $t=0$ based on equation~\eqref{eq:AB_vel_General} with the stored collision parameters $\sigma$. In this method, we also need to store all of the same information about colliding particles in the forward DSMC simulation as in the adjoint DSMC method, adding extra $28 N\mu T$ bytes needed for storage total of $8N + 28 N\mu T$ bytes.

(iv) For the direct discretization of the integrals in the adjoint equation~\eqref{backwards1}, we perform the forward simulation with the standard DSMC method. Afterwards, we need to convert velocity samples $\{v_{k,i}\}_{i=1}^N$ to a distribution function $f(v,t_k)$ via a histogram at every time step $t_k$ during the forward propagation. With $n_\text{grid}$ being the number of the grid points in each direction of the $v$-space, we have in total $n_\text{grid}^3$ grid points, which requires $8 n_\text{grid}^3 M$ bytes of memory to store the distribution function $f(v,t)$ at every time step. In our simulations, $n_\text{grid}=40$ considering the high computational complexity. With $n_\text{grid}^3=64000 \approx N/1000$ for $N=10^8$, we use $\approx1000$ particles per $\Delta v^3$ box on average.
During the backward solve, we need to store $n_\text{grid}^3$ values of the function $\gamma(v,t)$. For the right-hand side of~\eqref{backwards1_simplified2_discrete}, we can subtract one term of the right-hand side at a time to save memory. Therefore, in total, we need $8 n_\text{grid}^3 (M+2)$ bytes for the backward solve.
With $n_\text{grid}^3=N/1000$ and $M=20$, we get $8 n_\text{grid}^3 (M+2) = 0.176 N$ bytes  which is much smaller than the number of bytes required in the forward DSMC simulation.

\subsection{Error scaling}
To compare the accuracy, we assume that only one simulation is done for each value of the parameter $\alpha$. We use the forward Euler time integration for all the methods, so the time-integration error is $\cO(\Delta t)$ for all the four approaches. As it is hard to find all the constants in the error estimates, we only provide error scaling here.

(i) Due to the Monte Carlo representation of the distribution function $f(v,t)$, we have a spatial error $\cO(1/\sqrt{N})$ in computing the moments. For the gradient calculation, due to the finite difference scheme, we have errors $\cO((\Delta \alpha)^2)+\cO(1/(\Delta \alpha \sqrt{N}))$ contributed from~\eqref{eq:FD_error} and~\eqref{eq:rand_error}. 

(ii) Similar to the finite difference scheme, due to the Monte Carlo representation of the distribution function $f(v,t)$, we have a spatial error $\cO(1/\sqrt{N})$.

(iii) Again, the Monte Carlo representations of $f(v,t)$ and $\gamma_i(t_k)$ contribute to spatial error $\cO(1/\sqrt{N})$. The interpolation part of the DSMC-type scheme for the adjoint equation does not introduce significant extra error since the local error of a linear interpolation scales like $\propto (\overline{\Delta v})^2 \propto 1/N^{2/3} < 1/N^{1/2}$, where $\overline{\Delta v}$ is a characteristic distance between particles.    

(iv) In~\eqref{backwards1}, functions under the integral sign are smooth, non-singular (including the constant kernel $q$) and thus can be considered periodic in the $\varphi$-space, $\cos\theta$-space and the $v$-space (up to numerical round-off errors on the boundaries of the domain), so we can achieve effectively exponential convergence from the integral itself. However, since we have to compute $\gamma(v',t_{k+1})$ in \eqref{backwards1_simplified2_discrete} using interpolation, the overall order of convergence is limited by the order of interpolation. The linear interpolation introduces an error of size $\cO((\Delta v)^2)$. The forward simulation for~\eqref{eq:homoBoltz} can be done by a deterministic method by computing the integrals in a similar way here as to how we handle the continuous adjoint equation, which can then achieve a similar error of $\cO((\Delta v)^2)$. However, in this paper, we stick with computing $f(v,t)$ by simply converting the empirical distribution represented by the $N$ particles to a histogram at each time step. Thus, we get an additional error of $\cO(1/\sqrt{N})$ by the forward DSMC.

\subsection{Operation count}
In this subsection, we provide operation counts to illustrate the performance of the methods. The same as before, we only provide the dominant scaling for the operation counts.

(i) Colliding two particles takes 47 operations per collision pair in our code or 23.5 operations per colliding particle. That is, $23.5N_c M = 23.5 N\mu T$ total operations per simulation. To compute the gradient via the finite difference scheme, we need to run at least $D_\alpha+1$ simulations with slightly different values of the parameter $\alpha$, where $D_\alpha$ is the dimensionality of $\alpha$. Overall, to compute the gradient $\frac{d J_1}{d \alpha}$, we need $(D_\alpha+1) \times 23.5 N\mu T$ operations.

(ii) The adjoint DSMC algorithm takes 26 operations to collide two $\gamma$-particles together. There are $N_c/2$ pairs of particles collided at every time step, and $M$ total backward steps. We get the total operation count $26/2 N_c M = 13 N\mu T$. Also, $6$ more operations are needed to compute the post-collision unit direction $\sigma$ per a collision pair, and $3$ more operations are needed to compute the unit collision direction $(v-v_1){\hat ~}$ during the forward DSMC propagation for each collision pair. They are used later in the backward adjoint DSMC step. In total, we need $17.5 N\mu T$ operations per backward simulation in the adjoint DSMC algorithm.

(iii) In the DSMC-type scheme for the continuous adjoint equation~\eqref{backwards1}, we only need to perform $2$ operations and $1$ interpolation per particle at each time step. Interpolation is linear, but it is done here on a scattered (non-structured) set of data. We use a built-in MATLAB function to interpolate the scattered data, which uses a Delaunay triangulation on the scattered sample points to find neighboring points to the target location and perform the interpolation. The operation count of Delaunay triangulation scales like $\cO(N\log(N))$. This is why this method is relatively slow compared to the adjoint DSMC scheme. For the interpolation itself, we believe it is reasonable to assume that we need about 20 operations per point. As a result, it brings us roughly up to $(20 N_c + C N \log(N)) M = (20 + C \log(N)/\Delta t) N\mu T$ operations per simulation.
Besides, in the DSMC-type scheme, we need velocities $v_{k,i}$ at each time step $t_k$. If we do not save them during the forward propagation, we have to recover the velocity particles at a particular time from the final velocities $\{v_{F,i}\}_{i=1}^N$ by marching backward in time. This is done in the same way as how we backpropagate the adjoint particle $\bgamma_{k,i}$ from $t_{k+1}$ to $t_k$; see~\eqref{eq:AB_vel_General}. The back-propagation of the velocity particles requires an extra $17.5 N\mu T$ operations per simulation. The total number of operations is roughly $(40 + C \log(N)/\Delta t) N\mu T$ for the DSMC-type scheme.

(iv) Computing an integral over $dv$ and $d\sigma$ is the most expensive part at each time step for the direct discretization of~\eqref{backwards1}. Inside the integral, we need to compute $\gamma(v')$ using interpolation from the grid points in the $v$-space to the point $v'$. Again, we assume 20 operations are needed per grid point for the interpolation and the computation of $v'$.
The total number of operations needed to compute $\gamma(v,t_k)$ from $\gamma(v,t_{k+1})$ in \eqref{backwards1_simplified2_discrete} per grid point in the $v$-space is about $(21 n_{\varphi}n_{\theta}+4)n_\text{grid}^3$. Since there are $n_\text{grid}^3$ grid points in the $v$-space grid and $M$ time steps, the total operation count per simulation is about $(21 n_{\varphi}n_{\theta}+4)n_\text{grid}^6 M$.
Assuming $n_\text{grid}^3=N/1000$ and $n_{\varphi}=n_{\theta}=10$, we get $(21 n_{\varphi}n_{\theta}+4)n_\text{grid}^6 M = 0.0021 N^2 M$. There are drastically more operations in the direct discretization than the adjoint DSMC algorithm as a result of the $\cO(N^2)$ scaling where typically $N\geq10^6$.

\end{document}